\documentclass[11pt]{article}
\usepackage[left=1in,right=1in,top=1in,bottom=1in]{geometry}
\usepackage{times}
\usepackage{expl3}
\usepackage{cite}
\usepackage[table]{xcolor}
\usepackage{multirow}
\usepackage{stackengine} 
\usepackage{hhline}
\usepackage{lipsum}
\usepackage{titlesec}
\usepackage{wrapfig}
\usepackage{epsfig}
\usepackage{graphicx}
\usepackage{amsmath}
\usepackage[title]{appendix}
\usepackage{amssymb}
\usepackage{epstopdf}
\usepackage{boldline}
\usepackage{calligra}
\usepackage{url}
\usepackage{blindtext}

\newcommand{\define}{\stackrel{\mbox{\tiny def}}{=}}

\newtheorem{definition}{Definition}
\newtheorem{theorem}{Theorem}
\newtheorem{corollary}{Corollary}
\newtheorem{lemma}{Lemma}

\usepackage{mathtools}
\usepackage{epstopdf}
\usepackage{balance}
\usepackage{thmtools}
\usepackage{thm-restate}
\usepackage{hyperref}
\usepackage{cleveref}

\usepackage[ruled,vlined]{algorithm2e}
\include{pythonlisting}
\newcommand{\ostar}{\mathbin{\mathpalette\make@circled\star}}

\makeatletter
\newcommand{\removelatexerror}{\let\@latex@error\@gobble}
\makeatother
\makeatletter
\newcommand*{\rom}[1]{\expandafter\@slowromancap\romannumeral #1@}
\makeatother

\ExplSyntaxOn
\newcommand\latinabbrev[1]{
  \peek_meaning:NTF . {% Same as \@ifnextchar
    #1\@}%
  { \peek_catcode:NTF a {% Check whether next char has same catcode as \'a, i.e., is a letter
      #1.\@ }%
    {#1.\@}}}
\ExplSyntaxOff

\titleclass{\subsubsubsection}{straight}[\subsubsection]

\begin{document}
\vspace{1cm}
\title{Random Parametrization Double Tensors Integrals and Their Applications}\vspace{1.8cm}
\author{Shih~Yu~Chang 
% <-this % stops a space
\thanks{Shih Yu Chang is with the Department of Applied Data Science,
San Jose State University, San Jose, CA, U. S. A. (e-mail: {\tt
shihyu.chang@sjsu.edu}).
           }}

\maketitle

\begin{abstract}
In this work, we extend double tensor integrals (DTI) from our previous work to parametrization double tensors integrals (PDTI) by applying integral kernel transform bounds to upper bound PDTI norm and establishing a new perturbation formula. Besides, the convergence property of random PDTI is investigated and this property is utilized to characterize the relation between the original derivative tensor and the action result of PDTI to the original derivative tensor. These tools help us to derive new tail bounds for random tensors according to more general operator inequalities, e.g., Heinz inequality and Birman-Koplienko-Solomyak inequality. Moreover, new tail bounds about random tensors are also obtained according to our new derived perturbation formula and integral kernel transform bounds. 
% Birman-Koplienko-Solomyak
\end{abstract}


\section{Introduction}\label{sec:Introduction} 

% DTI

In order to consider the random tensor mean problem, we defined the notion about double tensor integrals (DTI) and discussed perturbation formula, Lipschitz estimation, and continuity issues for random DTI in~\cite{chang2022random}. Motivated by works in~\cite{azamov2009operator,potapov2009unbounded,potapov2010double} about applying double operator integration theory to noncommutative geometry, we extend DTI definition discussed in~\cite{chang2022random} to \emph{parametrization double tensors integrals} (PDTI). The idea to apply double operator integration techniques in the general area of operator inequalities can be traced back to the 1970s. For example, the original proof of Birman-Kopilenko-Solomyak inequality given in~\cite{birman1975estimates} depends on profound facts from double operator integration theory. The works from~\cite{potapov2009unbounded, potapov2010double} provide a framework by combining parametrization double operator integrals with Fourier transform bounds of perturbation function to prove various operator inequalities, e.g., Heinz inequality, Birman-Koplienko-Solomyak inequality, in a systematic approach.

In this work, we apply the framework from~\cite{potapov2009unbounded, potapov2010double} to random DTI. First, we extend operators from matrices format to tensors format by defining PDTI and consider more general integral kernel transform bounds, which will be used to upper bound PDTI norm. This will help us to associate the underlying perturbation function properties with PDTI norm estimation. Only Fourier transform is considered in~\cite{potapov2009unbounded, potapov2010double}. Second, we derive a more general perturbation formula, compared to Lemma 4 in~\cite{potapov2010double}, in Theorem~\ref{thm:perturbation theory}. Third, the convergence of random PDTI is provided by Lemma~\ref{lma:conv of T psi t}, which is used with Theorem~\ref{thm:perturbation theory} to characterize the relation between the original derivative tensor and the action result of PDTI to the original derivative tensor, see Lemma~\ref{lma:derivative after T psi action}. All these tools will help us to derive various new inequalities about random tensors. Originally, the Heinz inequality was proved in~\cite{heinz1951beitrage, mcintosh1979heinz, kosaki2011positive}. We extend the Heinz inequality by a tail bound format of random tensors in Theorem~\ref{thm:11lma}. Birman-Koplienko-Solomyak inequality was first proved in~\cite{birman1975estimates} with an alternative proof provided in~\cite{ando1988comparison}. Ando’s proof was later extended to semifinite von Neumann algebras in~\cite{dodds1995submajorization}. We extend this Birman-Kopilenko-Solomyak inequality to a more general setting by tail bounds of random tensors in Theorem~\ref{thm:16thm}. Other new inequalities are also obtained based on our new derived perturbation formula and integral kernel transform bounds, for example, Theorem~\ref{thm:17lma} and its corollary.

The rest of this paper is organized as follows. The terminologies related to tensors and fundamental facts about tensors are introduced in Section~\ref{sec:Fundamental of Tensors}. The extension of double tensor integrals, Parametrization Double Tensor Integrals (PDTI), is presented in Section~\ref{sec:Parametrizeation Double Tensor Integrals}. A new perturbation formula for a more general divided difference form is derived in Section~\ref{sec:Perturbation Formula}. In Section~\ref{sec:Limiting Behavior of Random Parametrizeation Double Tensor Integrals}, we will establish continuity conditions for PDTI using the convergence in mean for random tensors. In Section~\ref{sec:New Inequalities By PDTI}, we will apply the proposed PDTI to build several new inequalities of random tensors. Finally, conclusions will be drawn in Section~\ref{sec:Conclusions}.

\section{Fundamental of Tensors}\label{sec:Fundamental of Tensors}

Without loss of generality, one can partition the dimensions of a tensor into two groups, say $M$ and $N$ dimensions, separately. Thus, for two order-($M$+$N$) tensors: $\mathcal{X} \define (x_{i_1, \cdots, i_M, j_1, \cdots,j_N}) \in \mathbb{C}^{I_1 \times \cdots \times I_M\times
J_1 \times \cdots \times J_N}$ and $\mathcal{Y} \define (y_{i_1, \cdots, i_M, j_1, \cdots,j_N}) \in \mathbb{C}^{I_1 \times \cdots \times I_M\times
J_1 \times \cdots \times J_N}$, according to~\cite{chang2022random, chang2021convenient}, the \emph{tensor addition} $\mathcal{X} + \mathcal{Y}\in \mathbb{C}^{I_1 \times \cdots \times I_M\times
J_1 \times \cdots \times J_N}$ is given by 
\begin{eqnarray}\label{eq: tensor addition definition}
(\mathcal{X} + \mathcal{Y} )_{i_1, \cdots, i_M, j_1 , \cdots , j_N} &\define&
x_{i_1, \cdots, i_M, j_1 ,\cdots , j_N} \nonumber \\
& &+ y_{i_1, \cdots, i_M, j_1 ,\cdots , j_N}. 
\end{eqnarray}
On the other hand, for tensors $\mathcal{X} \define (x_{i_1, \cdots, i_M, j_1, \cdots,j_N}) \in \mathbb{C}^{I_1 \times \cdots \times I_M\times
J_1 \times \cdots \times J_N}$ and $\mathcal{Y} \define (y_{j_1, \cdots, j_N, k_1, \cdots,k_L}) \in \mathbb{C}^{J_1 \times \cdots \times J_N\times K_1 \times \cdots \times K_L}$, according to~\cite{chang2022random, chang2021convenient}, the \emph{Einstein product} (or simply referred to as \emph{tensor product} in this work) $\mathcal{X} \star_{N} \mathcal{Y} \in  \mathbb{C}^{I_1 \times \cdots \times I_M\times
K_1 \times \cdots \times K_L}$ is given by 
\begin{eqnarray}\label{eq: Einstein product definition}
\lefteqn{(\mathcal{X} \star_{N} \mathcal{Y} )_{i_1, \cdots, i_M,k_1, \cdots , k_L} \define} \nonumber \\ &&\sum\limits_{j_1, \cdots, j_N} x_{i_1, \cdots, i_M, j_1, \cdots,j_N}y_{j_1, \cdots, j_N, k_1, \cdots,k_L}. 
\end{eqnarray}

One can find more preliminary facts about tensors based on Einstein product in~\cite{chang2021convenient}. In the remaining of this paper, we will represent the scalar value $I_1 \times \cdots \times I_N$ by $\mathbb{I}_{1}^{N}$.

We also list other crucial tensor operations here. The \emph{trace} of a square tensor is equivalent to the summation of all diagonal entries such that 
\begin{eqnarray}\label{eq: tensor trace def}
\mathrm{Tr}(\mathcal{X}) \define \sum\limits_{1 \leq i_j \leq I_j,\hspace{0.05cm}j \in [M]} \mathcal{X}_{i_1, \dots, i_M,i_1, \dots, i_M}.
\end{eqnarray}
The \emph{inner product} of two tensors $\mathcal{X}$, $\mathcal{Y} \in \mathbb{C}^{I_1 \times \dots \times I_M\times J_1 \times \dots \times J_N}$ is given by 
\begin{eqnarray}\label{eq: tensor inner product def}
\langle \mathcal{X}, \mathcal{Y} \rangle \define \mathrm{Tr}\left(\mathcal{X}^H \star_M \mathcal{Y}\right).
\end{eqnarray}

From Theorem 3.2 in~\cite{liang2019further}, every Hermitian tensor $\mathcal{H} \in  \mathbb{C}^{I_1 \times \cdots \times I_N \times I_1 \times \cdots \times I_N}$ has the following decomposition
\begin{eqnarray}\label{eq:Hermitian Eigen Decom}
\mathcal{H} &=& \sum\limits_{i=1}^{\mathbb{I}_{1}^{N}} \lambda_i \mathcal{U}_i \star_1 \mathcal{U}^{H}_i  \mbox{
~with~~$\langle \mathcal{U}_i, \mathcal{U}_i \rangle =1$ and $\langle \mathcal{U}_i, \mathcal{U}_j \rangle = 0$ for $i \neq j$,} \nonumber \\
&\define& \sum\limits_{i=1}^{\mathbb{I}_{1}^{N}} \lambda_i \mathcal{P}_{\mathcal{U}_i}
\end{eqnarray}
where $ \mathcal{U}_i \in  \mathbb{C}^{I_1 \times \cdots \times I_N \times 1}$, and the tensor $\mathcal{P}_{\mathcal{U}_i}$ is defined as $\mathcal{U}_i \star_1 \mathcal{U}^{H}_i$. The values $\lambda_i$ are named as \emph{eigevalues}. A Hermitian tensor with the decomposition shown by Eq.~\eqref{eq:Hermitian Eigen Decom} is named as \emph{eigen-decomposition}. A Hermitian tensor $\mathcal{H}$ is a positive definite (or positive semi-definite) tensor if all its eigenvalues are positive (or nonnegative).  

\section{Parametrization Double Tensor Integrals}\label{sec:Parametrizeation Double Tensor Integrals}

Let $\psi: \mathbb{R} \times \mathbb{R} \rightarrow \mathbb{C}$ be a function with the following decomposition format in integrand as:
\begin{eqnarray}\label{eq1:def psi func}
\psi(\lambda_{\mathcal{A}}, \lambda_{\mathcal{B}}) = \int_{\Sigma} f_{\mathcal{A}, \sigma} (  \lambda_{\mathcal{A}} )  f_{\mathcal{B}, \sigma} (  \lambda_{\mathcal{B}} ) d \mu(\sigma),
\end{eqnarray}
where $\mu (\sigma)$ is a measure on measurable space $(\Sigma, \mu)$. Functions $ f_{\mathcal{A}, \sigma}: \mathbb{R} \rightarrow \mathbb{C}$ and $ f_{\mathcal{B}, \sigma}: \mathbb{R} \rightarrow \mathbb{C}$ are two bounded complex-valued functions satisfying 
\begin{eqnarray}\label{eq2:def psi func}
\int_{\Sigma} \left\Vert f_{\mathcal{A}, \sigma} (  \lambda_{\mathcal{A}} ) \right\Vert_{\infty} \left\Vert f_{\mathcal{B}, \sigma} (  \lambda_{\mathcal{B}} )  \right\Vert_{\infty} d \mu(\sigma) < \infty.
\end{eqnarray}
Let us collect all $\psi$ functions having the form as shown by Eq.~\eqref{eq1:def psi func} by a set $\Psi$ such that, for any given two functions $\psi_1, \psi_2 \in \Psi$ with   
\begin{eqnarray}
\psi_1(\lambda_{\mathcal{A}}, \lambda_{\mathcal{B}}) &=& \int_{\Sigma_1} f_{\mathcal{A}_1, \sigma_1} (  \lambda_{\mathcal{A}} )  f_{\mathcal{B}_1, \sigma_1} (  \lambda_{\mathcal{B}} ) d \mu_1(\sigma_1);\nonumber \\
\psi_2(\lambda_{\mathcal{A}}, \lambda_{\mathcal{B}}) &=& \int_{\Sigma_2} f_{\mathcal{A}_2, \sigma_2} (  \lambda_{\mathcal{A}} )  f_{\mathcal{B}_2, \sigma_2} (  \lambda_{\mathcal{B}} ) d \mu_2(\sigma_2),
\end{eqnarray}
we have new measure $(\Sigma_3, \mu_3)$ and new functions $ f_{\mathcal{A}_3, \sigma}, f_{\mathcal{B}_3, \sigma}$ satisfying Eq.~\eqref{eq2:def psi func} such that the following relation is valid~\footnote{In~\cite{potapov2010double}, this condition should be added 
to prove the Banach space of $\left\Vert \psi \right\Vert_{\Psi}$.} 
\begin{eqnarray}\label{eq3:def psi func}
 \int_{\Sigma_1} f_{\mathcal{A}_1, \sigma_1} (  \lambda_{\mathcal{A}} )  f_{\mathcal{B}_1, \sigma_1} (  \lambda_{\mathcal{B}} ) d \mu_1(\sigma_1) + \int_{\Sigma_2} f_{\mathcal{A}_2, \sigma_2} (  \lambda_{\mathcal{A}} )  f_{\mathcal{B}_2, \sigma_2} (  \lambda_{\mathcal{B}} ) d \mu_2(\sigma_2)    = \nonumber \\
\int_{\Sigma_3} f_{\mathcal{A}_3, \sigma_3} (  \lambda_{\mathcal{A}} )  f_{\mathcal{B}_3, \sigma_3} (  \lambda_{\mathcal{B}} ) d \mu_3(\sigma_3). 
\end{eqnarray}

We define the following norm function over the set $\Psi$ as
\begin{eqnarray}\label{eq1:def norm}
\left\Vert \psi \right\Vert_{\Psi} \define \min  \int_{\Sigma} \left\Vert f_{\mathcal{A}, \sigma} (  \lambda_{\mathcal{A}} ) \right\Vert_{\infty} \left\Vert   f_{\mathcal{B}, \sigma} (  \lambda_{\mathcal{B}} ) \right\Vert_{\infty}  d \mu(\sigma),
\end{eqnarray}
where the minimum is taken over all possible representations of Eq.~\eqref{eq1:def psi func}. With the condition provided by Eq.~\eqref{eq3:def psi func}, it is easy to verify that the norm defined by Eq.~\eqref{eq1:def norm} over the space $\Phi$ has the triangle inequality:
\begin{eqnarray}
\left\Vert \psi_1 +  \psi_2 \right\Vert_{\Psi} \leq \left\Vert \psi_1 \right\Vert_{\Psi} + \left\Vert \psi_2  \right\Vert_{\Psi}. 
\end{eqnarray}

Let $\mathcal{A}, \mathcal{B} \in \mathbb{C}^{I_1 \times \cdots \times I_N \times I_1 \times \cdots \times I_N} $ be Hermitian tensors with the following eigen-decompositions:
\begin{eqnarray}
\mathcal{A} &=& \sum\limits_{i=1}^{\mathbb{I}_{1}^{N}} \lambda_{\mathcal{A}, i} \mathcal{U}_{\mathcal{A}, i} \star_1 \mathcal{U}^{H}_{\mathcal{A}, i}
\define \sum\limits_{i=1}^{\mathbb{I}_{1}^{N}} \lambda_{\mathcal{A}, i} \mathcal{P}_{\mathcal{A}, i},
\end{eqnarray}
and
\begin{eqnarray}
\mathcal{B} &=& \sum\limits_{j=1}^{\mathbb{I}_{1}^{N}} \lambda_{\mathcal{B}, j} \mathcal{U}_{\mathcal{B}, j} \star_1 \mathcal{U}^{H}_{\mathcal{B}, j}
\define \sum\limits_{j=1}^{\mathbb{I}_{1}^{N}} \lambda_{\mathcal{B}, j} \mathcal{P}_{\mathcal{B}, j},
\end{eqnarray}
where $\mathcal{P}_{\mathcal{A}, i}$ and $\mathcal{P}_{\mathcal{B}, j}$ are projection tensors of tensors $\mathcal{A}$ and $\mathcal{B}$, respectively. We also have the function $\psi(\lambda_{\mathcal{A}}, \lambda_{\mathcal{B}})$ associated to eigenvalues of $\lambda_{\mathcal{A}}$ and $\lambda_{\mathcal{B}}$ defined by Eq.~\eqref{eq1:def psi func}. Then, we can define a parametrize double tensor integrals (PDTI) over the measurable space $(\Sigma, \mu)$, represented by $T_{\psi}(\mathcal{X})$, as:
\begin{eqnarray}\label{eq1:def pDTI}
T_{\psi}(\mathcal{X}) &=& \int_{\Sigma}  \left(\sum\limits_{i=1}^{\mathbb{I}_1^N}  f_{\mathcal{A}, \sigma} \left(  \lambda_{\mathcal{A}, i} \right)  \mathcal{P}_{\mathcal{A}, i} \right)\star_N \mathcal{X} \star_N \left(\sum\limits_{j=1}^{\mathbb{I}_1^N}  f_{\mathcal{B}, \sigma} \left(  \lambda_{\mathcal{B}, j} \right) \mathcal{P}_{\mathcal{B}, j} \right)  d \mu(\sigma)
\end{eqnarray}
where $\mathcal{X} \in \mathbb{C}^{I_1 \times \cdots \times I_N \times I_1 \times \cdots \times I_N}$. $T_{\psi}(\mathcal{X})$ is called a random PDTI if $\lambda_{\mathcal{A}, i}, \lambda_{\mathcal{B}}$ are random variables and $\mathcal{P}_{\mathcal{A}, i}, \mathcal{P}_{\mathcal{B}, j}$ are random tensors.

From the definition provided by Eq.~\eqref{eq1:def pDTI}, we have the following Lemma about $T_{\psi}(\mathcal{X}) $. 

\begin{lemma}[Kernel of the mapping $\psi \rightarrow T_{\psi}$ is zero]\label{lma:kernel zero}
Given the function $\psi(\lambda_{\mathcal{A}}, \lambda_{\mathcal{B}})$ defined by Eq.~\eqref{eq1:def psi func}, the Kernel space of the mapping $\psi \rightarrow T_{\psi}$ is zero.
\end{lemma}
\textbf{Proof:}
It is enough to prove that if functions $ f_{\mathcal{A}, \sigma}$ and $ f_{\mathcal{B}, \sigma}$ have the following property:
\begin{eqnarray}\label{eq1:lma:kernel zero}
\int_{\Sigma} f_{\mathcal{A}, \sigma} (  \lambda_{\mathcal{A}} )  f_{\mathcal{B}, \sigma} (  \lambda_{\mathcal{B}} ) d \mu(\sigma) = 0,
\end{eqnarray}
we have
\begin{eqnarray}\label{eq2:lma:kernel zero}
\mathrm{Tr}\left(T_{\psi}( \mathcal{X}) \star_N \mathcal{Y}  \right) = 0,
\end{eqnarray}
where $\mathcal{X}$ and $\mathcal{Y}$ are any tensors with dimensions $\mathbb{C}^{I_1 \times \cdots \times I_N \times I_1 \times \cdots \times I_N}$.  

Suppose we have the following expression for tensors $\mathcal{X}$ and $\mathcal{Y}$:
\begin{eqnarray}
\mathcal{X} = \mathcal{U}_{A} \star_1 \mathcal{V}_{B},
\end{eqnarray}
where $\mathcal{U}_{A}, \mathcal{V}_{B} \in \mathbb{C}^{I_1 \times \cdots \times I_N}$; and
\begin{eqnarray}
\mathcal{Y} = \mathcal{U}_{B} \star_1 \mathcal{V}_{A},
\end{eqnarray}
where $\mathcal{U}_{B}, \mathcal{V}_{A} \in \mathbb{C}^{I_1 \times \cdots \times I_N}$. For any $\sigma \in \Sigma$, we have
\begin{eqnarray}\label{eq3:lma:kernel zero}
\lefteqn{\mathrm{Tr}\left( \left(\sum\limits_{i=1}^{\mathbb{I}_1^N}  f_{\mathcal{A}, \sigma} \left(  \lambda_{\mathcal{A}, i} \right)  \mathcal{P}_{\mathcal{A}, i} \right)\star_N \mathcal{X} \star_N \left(\sum\limits_{j=1}^{\mathbb{I}_1^N}  f_{\mathcal{B}, \sigma} \left(  \lambda_{\mathcal{B}, j} \right) \mathcal{P}_{\mathcal{B}, j} \right) \star_N \mathcal{Y} \right) }\nonumber \\
&& = \left\langle \sum\limits_{i=1}^{\mathbb{I}_1^N}  f_{\mathcal{A}, \sigma} \left(  \lambda_{\mathcal{A}, i} \right)  \mathcal{P}_{\mathcal{A}, i}  \star_N \mathcal{U}_{A} ,  \mathcal{V}_{A} \right\rangle 
\left\langle \sum\limits_{j=1}^{\mathbb{I}_1^N}  f_{\mathcal{B}, \sigma} \left(  \lambda_{\mathcal{B}, j} \right)  \mathcal{P}_{\mathcal{B}, j}  \star_N \mathcal{U}_{B} ,  \mathcal{V}_{B} \right\rangle  \nonumber \\ 
&& = \sum\limits_{i=1}^{\mathbb{I}_1^N} \sum\limits_{j=1}^{\mathbb{I}_1^N}
\left( f_{\mathcal{A}, \sigma}  \left(  \lambda_{\mathcal{A}, i} \right)  f_{\mathcal{B}, \sigma}  \left(  \lambda_{\mathcal{B}, j} \right) \right) \left( \left\langle \mathcal{P}_{\mathcal{A}, i}  \star_N \mathcal{U}_{A} ,  \mathcal{V}_{A} \right\rangle  
\left\langle \mathcal{P}_{\mathcal{B}, j}  \star_N \mathcal{U}_{B} ,  \mathcal{V}_{B} \right\rangle     \right).
\end{eqnarray}

If we integrate both sides at Eq.~\eqref{eq3:lma:kernel zero} with respect to $\sigma$, we have 
\begin{eqnarray}\label{eq4:lma:kernel zero}
\mathrm{Tr}\left(T_{\psi}\left( \mathcal{X}\right) \star_N \mathcal{Y} \right) &=& 
\sum\limits_{i=1}^{\mathbb{I}_1^N} \sum\limits_{j=1}^{\mathbb{I}_1^N}
\left( f_{\mathcal{A}, \sigma}  \left(  \lambda_{\mathcal{A}, i} \right)  f_{\mathcal{B}, \sigma}  \left(  \lambda_{\mathcal{B}, j} \right) \right) 
\overbracket{\left[\int_{\Sigma}   \left( f_{\mathcal{A}, \sigma}  \left(  \lambda_{\mathcal{A}, i} \right)  f_{\mathcal{B}, \sigma}  \left(  \lambda_{\mathcal{B}, j} \right) \right) d \mu (\sigma) \right]}^{\psi (    \lambda_{\mathcal{A}, i}, \lambda_{\mathcal{B}, j} )} \nonumber \\
&& 
\left( \left\langle \mathcal{P}_{\mathcal{A}, i}  \star_N \mathcal{U}_{A} ,  \mathcal{V}_{A} \right\rangle  
\left\langle \mathcal{P}_{\mathcal{B}, j}  \star_N \mathcal{U}_{B} ,  \mathcal{V}_{B} \right\rangle     \right).
\end{eqnarray}
Then, if the function $\psi$ becomes  $0$, we have $\mathrm{Tr}\left(T_{\psi}\left( \mathcal{X}\right) \star_N \mathcal{Y} \right) = 0$. This indicates that $T_{\psi}$ will be zero. 
$\hfill \Box$

Our next lemma is about the norm estimate of $T_{\psi}$. The spectral norm of a tensor is assumed here, i.e., $\left\Vert \mathcal{X} \right\Vert =  s_{\max }(\mathcal{A})$, where $s_{\max}$ represents the largest singular value of the tensor $\mathcal{A}$, see Theorem 3.2 in~\cite{liang2019further} about the singular values definition of a tensor.

\begin{lemma}[Norm estimate of $T_{\psi}$ by $\psi$ norm]\label{lma:T psi by psi}
Let $T_{\psi}(\mathcal{X})$ defined by Eq.~\eqref{eq1:def pDTI}, we have the following spectral norm estimate 
\begin{eqnarray}
\left\Vert T_{\psi}(\mathcal{X}) \right\Vert \leq \left(\mathbb{I}_1^N \right)^2 \left\Vert \psi \right\Vert_{\Psi} \left\Vert \mathcal{X} \right\Vert.
\end{eqnarray}
\end{lemma}
\textbf{Proof:}
Suppose we select a  $\psi \in \Psi$ and $\epsilon > 0$ such that 
\begin{eqnarray}
\int_{\Sigma} \left\Vert f_{\mathcal{A}, \sigma} (  \lambda_{\mathcal{A}} ) \right\Vert_{\infty} \left\Vert   f_{\mathcal{B}, \sigma} (  \lambda_{\mathcal{B}} ) \right\Vert_{\infty}  d \mu(\sigma) < 
 \left( \left\Vert \psi \right\Vert_{\Psi} + \epsilon \right).
\end{eqnarray} 

We also have
\begin{eqnarray}
\lefteqn{\left\Vert \left(\sum\limits_{i=1}^{\mathbb{I}_1^N}  f_{\mathcal{A}, \sigma} \left(  \lambda_{\mathcal{A}, i} \right)  \mathcal{P}_{\mathcal{A}, i} \right)\star_N \mathcal{X} \star_N \left(\sum\limits_{j=1}^{\mathbb{I}_1^N}  f_{\mathcal{B}, \sigma} \left(  \lambda_{\mathcal{B}, j} \right) \mathcal{P}_{\mathcal{B}, j} \right)   \right\Vert } \nonumber \\
&&  \leq_1 \left\Vert \sum\limits_{i=1}^{\mathbb{I}_1^N}  f_{\mathcal{A}, \sigma} \left(  \lambda_{\mathcal{A}, i} \right)  \mathcal{P}_{\mathcal{A}, i} \right\Vert  \left\Vert \mathcal{X} \right\Vert   \left\Vert \sum\limits_{j=1}^{\mathbb{I}_1^N}  f_{\mathcal{B}, \sigma} \left(  \lambda_{\mathcal{B}, j} \right)  \mathcal{P}_{\mathcal{B}, j} \right\Vert \nonumber \\
&&  \leq_2  \left(\mathbb{I}_1^N \right)^2  \left\Vert f_{\mathcal{A}, \sigma} (  \lambda_{\mathcal{A}} ) \right\Vert_{\infty} \left\Vert   f_{\mathcal{B}, \sigma} (  \lambda_{\mathcal{B}} ) \right\Vert_{\infty}  \left\Vert \mathcal{X} \right \Vert,
\end{eqnarray} 
where the inequality $\leq_1$ is based on the submultiplicative of spectral norm and the inequality $\leq_2$ is based on the triangle inequality and the fact that the spectral norm of $\mathcal{P}_{\mathcal{A}, i}$ and $\mathcal{P}_{\mathcal{B}, j}$ are one.

Then, we can have the following relation
\begin{eqnarray}
\lefteqn{\left\Vert T_{\psi} (\mathcal{X}) \right\Vert \leq \int_{\Sigma} \left\Vert \left(\sum\limits_{i=1}^{\mathbb{I}_1^N}  f_{\mathcal{A}, \sigma} \left(  \lambda_{\mathcal{A}, i} \right)  \mathcal{P}_{\mathcal{A}, i} \right)\star_N \mathcal{X} \star_N \left(\sum\limits_{j=1}^{\mathbb{I}_1^N}  f_{\mathcal{B}, \sigma} \left(  \lambda_{\mathcal{B}, j} \right) \mathcal{P}_{\mathcal{B}, j} \right)   \right\Vert d \mu(\sigma)}  \nonumber \\
&&  \leq  \left(\mathbb{I}_1^N \right)^2  \left[ \int_{\Sigma} \left\Vert f_{\mathcal{A}, \sigma} (  \lambda_{\mathcal{A}} ) \right\Vert_{\infty} \left\Vert   f_{\mathcal{B}, \sigma} (  \lambda_{\mathcal{B}} ) \right\Vert_{\infty}     \right] \left\Vert \mathcal{X} \right\Vert  \nonumber \\
&&  \leq \left(\mathbb{I}_1^N \right)^2 \left(  \left\Vert \psi \right\Vert_{\Psi} + \epsilon \right) \left\Vert \mathcal{X} \right\Vert .
\end{eqnarray} 
This Lemma is proved by taking $\epsilon \rightarrow 0$. 
$\hfill \Box$

From Lemma~\ref{lma:T psi by psi}, we only bound the PDTI in terms of $\left\Vert \psi \right\Vert_{\Psi}$. Following theorem will give the bound for $\left\Vert \psi \right\Vert_{\Psi}$ by the property of $\psi$ function. 

\begin{theorem}\label{thm:bound for psi norm}
Suppose we are given an integral transform as:
\begin{eqnarray}\label{eq2:thm:bound for psi norm}
g(t) = \int_{\mathbb{R}} K(s,t) \tilde{g}(s)  d s.
\end{eqnarray}
If the variable $t$ is associated to eigenvalues of $\lambda_{\mathcal{A}}$ and $\lambda_{\mathcal{B}}$ by the following bivariable function as 
\begin{eqnarray}
t = \beta(\lambda_{\mathcal{A}}, \lambda_{\mathcal{B}}), 
\end{eqnarray}
and $\psi(\lambda_{\mathcal{A}}, \lambda_{\mathcal{B}})$ is assumed to be expressed as
\begin{eqnarray}
\psi(\lambda_{\mathcal{A}}, \lambda_{\mathcal{B}}) &=& g ( \beta(\lambda_{\mathcal{A}}, \lambda_{\mathcal{B}}) ) \nonumber \\
&=& \int_{\mathbb{R}} K(s,   \beta(\lambda_{\mathcal{A}}, \lambda_{\mathcal{B}})    ) \tilde{g}(s)  d s  \nonumber \\
&=& 
 \int_{\mathbb{R}} f_{\mathcal{A}, s }( \lambda_{\mathcal{A}} )    
 f_{\mathcal{B}, s }( \lambda_{\mathcal{B}} )   \tilde{g}(s)  d s,
\end{eqnarray}
where $ f_{\mathcal{A}, s }( \lambda_{\mathcal{A}} )  =  f_{\mathcal{A}, \sigma }( \lambda_{\mathcal{A}} ) $ and $ f_{\mathcal{B}, s }( \lambda_{\mathcal{B}} )  =  f_{\mathcal{B}, \sigma }( \lambda_{\mathcal{B}} ) $. For all $\sigma \in \Sigma$, we assume that $\left\Vert f_{\mathcal{A}, \sigma} \right\Vert_{\infty} \leq c_{\mathcal{A}}$ and $\left\Vert f_{\mathcal{B}, \sigma} \right\Vert_{\infty} \leq c_{\mathcal{B}}$, where both $c_{\mathcal{A}}$ and $c_{\mathcal{B}}$ are two positive real numbers.

Then, we have 
\begin{eqnarray}\label{eq1:thm:bound for psi norm}
\left\Vert \psi\right\Vert_{\Psi} &\leq&  c_{\mathcal{A}}c_{\mathcal{B}}  \left(  \int_{\mathbb{R}} \left( \max\limits_{t}\left\vert K(s,t) \right\vert \right)  d s  \right) \left\Vert g(t) \right\Vert_{\infty}.
\end{eqnarray}
\end{theorem}
\textbf{Proof:}
From the definition of $\left\Vert \psi\right\Vert_{\Psi}$, we have
\begin{eqnarray}
\left\Vert \psi\right\Vert_{\Psi} &\leq&  \int_{\Sigma} \left\Vert f_{\mathcal{A}, \sigma} (  \lambda_{\mathcal{A}} ) \right\Vert_{\infty} \left\Vert   f_{\mathcal{B}, \sigma} (  \lambda_{\mathcal{B}} ) \right\Vert_{\infty}  d \mu(\sigma)   \nonumber \\
&\leq_1& c_{\mathcal{A}}c_{\mathcal{B}}  \int_{\Sigma}   d \mu(\sigma) 
  \nonumber \\
&=_2&  c_{\mathcal{A}}c_{\mathcal{B}}  \int_{\mathbb{R}} \left\vert  \tilde{g}(s)  \right\vert  d s  
  \nonumber \\
&\leq_3&  c_{\mathcal{A}}c_{\mathcal{B}}  \left(  \int_{\mathbb{R}} \left( \max\limits_{t}\left\vert K(s,t) \right\vert \right)  d s  \right) \left\Vert g(t) \right\Vert_{\infty},
\end{eqnarray}
where the inequality $\leq_1$ comes from assumptions about $ \left\Vert f_{\mathcal{A}, \sigma} (  \lambda_{\mathcal{A}} ) \right\Vert_{\infty}$ and $ \left\Vert   f_{\mathcal{B}, \sigma} (  \lambda_{\mathcal{B}} ) \right\Vert_{\infty}$, the equality $=_2$ is obtained by setting $\Sigma = \mathbb{R}$ and $d \mu (\sigma) =  \left\vert  \tilde{g}(s)  \right\vert  d s  $, and the inequality $\leq_3$ comes from H\"{o}lder's inequality with $p = \infty, q=1$. This theorem is proved.
$\hfill \Box$

We will have following corollaries according to Theorem~\ref{thm:bound for psi norm} by 
choosing different transform functions $K(s,t)$. But, we need the following Lemma about the $L^1$ estimate of Fourier transform. 

\begin{lemma}\label{lma:FT bound}
If $g(t): \mathbb{R} \rightarrow \mathbb{C}$ is an an absolutely continuous function with $g, g'$ are $L^2$ function, we have
\begin{eqnarray}
\left\Vert \tilde{g}(s) \right\Vert_1 \leq \min\limits_{c > 0}\left(\sqrt{2c} \left\Vert g(t) \right\Vert_{2} + \sqrt{2/c} \left\Vert g'(t) \right\Vert_{2}
  \right), 
\end{eqnarray}
where $c$ is any positive real number and $\tilde{g}(s)$ is the Fourier transform of $g(t)$.
\end{lemma}
\textbf{Proof:}
Since we have 
\begin{eqnarray}
\int_{\mathbb{R}} \left\vert \tilde{g}(s) \right\vert ds &=& 
\int_{s \in [-c, c]} \left\vert \tilde{g}(s) \right\vert ds + 
\int_{s \notin [-c, c]}  \left\vert s \right\vert^{-1}  \left\vert s \tilde{g}(s) \right\vert ds \nonumber \\
&\leq_1 & \sqrt{2c}\left( \int_{s \in [-c, c]} \left\vert \tilde{g}(s) \right\vert^2 ds     \right)^{1/2}
\nonumber \\
& & + \left(  \int_{s \notin [-c, c]} \left\vert s \right\vert^{-2} ds      \right)^{1/2} \cdot  
 \left(  \int_{s \notin [-c, c]} \left\vert s \tilde{g}(s) \right\vert^{2} ds      \right)^{1/2} \nonumber \\
&\leq_2 & \sqrt{2c} \left\Vert g(t) \right\Vert_{2} + \sqrt{2/c} \left\Vert g'(t) \right\Vert_{2}
\end{eqnarray}
where $\leq_1$ comes from the Cauchy–Schwarz inequality, and $\leq_2$ uses Plancherel
identity and the $L^2$ norm has larger support $\mathbb{R}$ than $s \notin [-c, c]  $. This Lemma is proved by taking the minimization over the positive variable $c$. 
$\hfill \Box$

\begin{corollary}\label{cor:bound for psi norm Fourier}
Suppose we are given a Fourier transform 
\begin{eqnarray}\label{eq2:cor:bound for psi norm Fourier}
g(t) = \int_{\mathbb{R}} \tilde{g}(s)  e^{\iota t s}   d s,
\end{eqnarray}
where $\iota = \sqrt{-1}$.  If the variable $t$ is associated to eigenvalues of $\lambda_{\mathcal{A}}$ and $\lambda_{\mathcal{B}}$ by the following bivariable function as 
\begin{eqnarray}
t = \log \left( \frac{\gamma(\lambda_{\mathcal{A}}) }{ \kappa(\lambda_{\mathcal{B}})} \right), 
\end{eqnarray}
where $\gamma: \mathbb{R} \rightarrow \mathbb{R}^{+}$ and $ \kappa: \mathbb{R} \rightarrow \mathbb{R}^{+}$. 

If $\psi(\lambda_{\mathcal{A}}, \lambda_{\mathcal{B}}   ) = g\left( \log \left( \frac{\gamma(\lambda_{\mathcal{A}}) }{ \kappa(\lambda_{\mathcal{B}})} \right) \right)$, then, we have 
\begin{eqnarray}\label{eq1:thm:bound for psi norm Fourier}
\left\Vert \psi\right\Vert_{\Psi} &\leq&  \min\limits_{c > 0}\left(\sqrt{2c} \left\Vert g(t) \right\Vert_{2} + \sqrt{2/c} \left\Vert g'(t) \right\Vert_{2}
  \right), 
\end{eqnarray}
where $\left\Vert ~ \right\Vert_{2}$ is $L^2$ function norm. 
\end{corollary}
\textbf{Proof:}
Since $t = \log \left( \frac{\gamma(\lambda_{\mathcal{A}}) }{ \kappa(\lambda_{\mathcal{B}})} \right)$, we have 
\begin{eqnarray}
\psi(\lambda_{\mathcal{A}}, \lambda_{\mathcal{B}}   ) = g\left( \log \left( \frac{ \gamma(\lambda_{\mathcal{A}}) }{ \kappa(\lambda_{\mathcal{B}})} \right) \right)
= \int_{\mathbb{R}}  \tilde{g}(s) (  \gamma(\lambda_{\mathcal{A}})  )^{  \iota s }   (  \kappa(\lambda_{\mathcal{B}})  )^{ - \iota s }      d s. 
\end{eqnarray}

If we set the following parameters: $\Sigma = \mathbb{R}$, $d \mu (\sigma) = \left\vert \tilde{g}(s) \right\vert ds $, $   f_{\mathcal{A}, \sigma} (  \lambda_{\mathcal{A}} ) =  (  \gamma(\lambda_{\mathcal{A}})  )^{  \iota s }     $ and $f_{\mathcal{B}, \sigma} (  \lambda_{\mathcal{B}} ) = \kappa(\lambda_{\mathcal{B}})  )^{ - \iota s }   $, we obtain
\begin{eqnarray}
\left\Vert \psi\right\Vert_{\Psi} &\leq&  \int_{\Sigma} \left\Vert f_{\mathcal{A}, \sigma} (  \lambda_{\mathcal{A}} ) \right\Vert_{\infty} \left\Vert   f_{\mathcal{B}, \sigma} (  \lambda_{\mathcal{B}} ) \right\Vert_{\infty}  d \mu(\sigma)   \nonumber \\
&\leq_1& 1 \times  1 \times  \int_{\Sigma}   d \mu(\sigma) 
  \nonumber \\
&=&  1 \times  1 \times  \int_{\mathbb{R}} \left\vert  \tilde{g}(s)  \right\vert  d s  
  \nonumber \\
&\leq_2& \min\limits_{c > 0}\left(\sqrt{2c} \left\Vert g(t) \right\Vert_{2} + \sqrt{2/c} \left\Vert g'(t) \right\Vert_{2}
  \right), 
\end{eqnarray}
where the inequality $\leq_1$ comes from assumptions about $ \left\Vert f_{\mathcal{A}, \sigma} (  \lambda_{\mathcal{A}} ) \right\Vert_{\infty} = 1$ and $\leq_2$ comes from Lemma~\ref{lma:FT bound}. This Corollary is proved.
$\hfill \Box$

\begin{corollary}\label{thm:bound for psi norm Laplace}
Suppose we are given a transform 
\begin{eqnarray}\label{eq2:thm:bound for psi norm Laplace}
g(t) = \int_{\mathbb{R}} \tilde{g}(s)  e^{\alpha t +   \iota t s}   d s.
\end{eqnarray}
If the variable $t$ is associated to eigenvalues of $\lambda_{\mathcal{A}}$ and $\lambda_{\mathcal{B}}$ by the following bivariable function as 
\begin{eqnarray}
t = \log \left( \frac{\gamma(\lambda_{\mathcal{A}}) }{ \kappa(\lambda_{\mathcal{B}})} \right), 
\end{eqnarray}
where $\gamma: \mathbb{R} \rightarrow \mathbb{R}^{+}$ and $ \kappa: \mathbb{R} \rightarrow \mathbb{R}^{+}$. 

If $\psi(\lambda_{\mathcal{A}}, \lambda_{\mathcal{B}}   ) = g\left( \log \left( \frac{\gamma(\lambda_{\mathcal{A}}) }{ \kappa(\lambda_{\mathcal{B}})} \right) \right)$, then, we have 
\begin{eqnarray}\label{eq1:thm:bound for psi norm Laplace}
\left\Vert \psi\right\Vert_{\Psi} &\leq&   \gamma^{\alpha}(\lambda^*_{\mathcal{A}})   \kappa^{-\alpha}(\lambda^*_{\mathcal{B}})   \min\limits_{c > 0}\left(\sqrt{2c} \left\Vert g(t) \right\Vert_{2} + \sqrt{2/c} \left\Vert g'(t) \right\Vert_{2}
  \right), 
\end{eqnarray}
where $\gamma^{\alpha}(\lambda^*_{\mathcal{A}})$ and $\kappa^{-\alpha}(\lambda^*_{\mathcal{B}})$ are the maximum values of the functions $\gamma^{\alpha}(\lambda_{\mathcal{A}})$ and $\kappa^{-\alpha}(\lambda_{\mathcal{B}})$, respectively. 
\end{corollary}
\textbf{Proof:}
Since $t = \log \left( \frac{\gamma(\lambda_{\mathcal{A}}) }{ \kappa(\lambda_{\mathcal{B}})} \right)$, we have 
\begin{eqnarray}
\psi(\lambda_{\mathcal{A}}, \lambda_{\mathcal{B}}   ) = g\left( \log \left( \frac{ \gamma(\lambda_{\mathcal{A}}) }{ \kappa(\lambda_{\mathcal{B}})} \right) \right)
= \int_{\mathbb{R}}  \tilde{g}(s) \gamma^{\alpha}(\lambda_{\mathcal{A}}) (  \gamma(\lambda_{\mathcal{A}})  )^{  \iota s }  \kappa^{-\alpha}(\lambda_{\mathcal{B}})  (  \kappa(\lambda_{\mathcal{B}})  )^{ - \iota s }      d s. 
\end{eqnarray}

If we set the following parameters: $\Sigma = \mathbb{R}$, $d \mu (\sigma) = \left\vert \tilde{g}(s) \right\vert ds $, $   f_{\mathcal{A}, \sigma} (  \lambda_{\mathcal{A}} ) = \gamma^{\alpha}(\lambda_{\mathcal{A}})   (  \gamma(\lambda_{\mathcal{A}})  )^{  \iota s }     $ and $f_{\mathcal{B}, \sigma} (  \lambda_{\mathcal{B}} ) =\kappa^{-\alpha}(\lambda_{\mathcal{B}})  (\kappa(\lambda_{\mathcal{B}})  )^{ - \iota s }   $, we obtain
\begin{eqnarray}
\left\Vert \psi\right\Vert_{\Psi} &\leq&  \int_{\Sigma} \left\Vert f_{\mathcal{A}, \sigma} (  \lambda_{\mathcal{A}} ) \right\Vert_{\infty} \left\Vert   f_{\mathcal{B}, \sigma} (  \lambda_{\mathcal{B}} ) \right\Vert_{\infty}  d \mu(\sigma)   \nonumber \\
&\leq_1& \gamma^{\alpha}(\lambda^*_{\mathcal{A}})   \kappa^{-\alpha}(\lambda^*_{\mathcal{B}})  \int_{\Sigma}   d \mu(\sigma) 
  \nonumber \\
&=&  \gamma^{\alpha}(\lambda^*_{\mathcal{A}})   \kappa^{-\alpha}(\lambda^*_{\mathcal{B}})  \int_{\mathbb{R}} \left\vert  \tilde{g}(s)  \right\vert  d s  
  \nonumber \\
&\leq_2& \gamma^{\alpha}(\lambda^*_{\mathcal{A}})   \kappa^{-\alpha}(\lambda^*_{\mathcal{B}})  \min\limits_{c > 0}\left(\sqrt{2c} \left\Vert g(t) \right\Vert_{2} + \sqrt{2/c} \left\Vert g'(t) \right\Vert_{2}
  \right), 
\end{eqnarray}
where the inequality $\leq_1$ comes from the definition of $\gamma^{\alpha}(\lambda^*_{\mathcal{A}})$ and $\kappa^{-\alpha}(\lambda^*_{\mathcal{B}})$, and $\leq_2$ comes from Lemma~\ref{lma:FT bound}. This Corollary is also proved.
$\hfill \Box$

\section{Perturbation Formula}\label{sec:Perturbation Formula} 

The main purpose of this section is to prepare a perturbation formula for the tensor operator $T_{\psi}$ with respect to a more general divided difference form. We begin with some preparation lemmas. 

\begin{lemma}\label{lma:homomorphism}
The mapping $\psi \rightarrow T_{\psi}$ is a homomorphism.
\end{lemma}
\textbf{Proof:}
We note that $\Psi$ is a Banach algebra since it is closed under the multiplication and it is also continuous with respect to the norm of $\Phi$ defined by Eq.~\eqref{eq1:def norm}. 

We define $\psi_1$ and $\psi_2$ as follows
\begin{eqnarray}
\psi_1(\lambda_{\mathcal{A}}, \lambda_{\mathcal{B}}) &=& \int_{\Sigma_1} f_{\mathcal{A}, \sigma_1} (  \lambda_{\mathcal{A}} )  f_{\mathcal{B}, \sigma_1} (  \lambda_{\mathcal{B}} ) d \mu_1(\sigma_1);\nonumber \\
\psi_2(\lambda_{\mathcal{A}}, \lambda_{\mathcal{B}}) &=& \int_{\Sigma_2} f_{\mathcal{A}, \sigma_2} (  \lambda_{\mathcal{A}} )  f_{\mathcal{B}, \sigma_2} (  \lambda_{\mathcal{B}} ) d \mu_2(\sigma_2),
\end{eqnarray}
and assume that $\psi_3$ is the product of $\psi_1$ and $\psi_2$. Then, we can further define the following terms:
\begin{eqnarray}
F_{\mathcal{A}, \sigma_1} &\define& \sum\limits_{i=1}^{\mathbb{I}_1^N}  f_{\mathcal{A}, \sigma_1} \left(  \lambda_{\mathcal{A}, i} \right)  \mathcal{P}_{\mathcal{A}, i}, ~~~
F_{\mathcal{B}, \sigma_1} \define \sum\limits_{i=1}^{\mathbb{I}_1^N}  f_{\mathcal{B}, \sigma_1} \left(  \lambda_{\mathcal{B}, i} \right)  \mathcal{P}_{\mathcal{B}, i}, \nonumber \\
F_{\mathcal{A}, \sigma_2} &\define& \sum\limits_{i=1}^{\mathbb{I}_1^N}  f_{\mathcal{A}, \sigma_2} \left(  \lambda_{\mathcal{A}, i} \right)  \mathcal{P}_{\mathcal{A}, i} ,~~~
F_{\mathcal{B}, \sigma_2} \define \sum\limits_{i=1}^{\mathbb{I}_1^N}  f_{\mathcal{B}, \sigma_2} \left(  \lambda_{\mathcal{B}, i} \right)  \mathcal{P}_{\mathcal{B}, i}.
\end{eqnarray}
From the spectral mapping theorem, we have
\begin{eqnarray}
F_{\mathcal{A}, \sigma_1} \star_N F_{\mathcal{A}, \sigma_2} &= & \sum\limits_{i=1}^{\mathbb{I}_1^N}  f_{\mathcal{A}, \sigma_1} \left(  \lambda_{\mathcal{A}, i} \right) f_{\mathcal{A}, \sigma_2} \left(  \lambda_{\mathcal{A}, i} \right)   \mathcal{P}_{\mathcal{A}, i}, \nonumber \\
F_{\mathcal{B}, \sigma_1} \star_N F_{\mathcal{B}, \sigma_2} &= & \sum\limits_{i=1}^{\mathbb{I}_1^N}  f_{\mathcal{B}, \sigma_1} \left(  \lambda_{\mathcal{B}, i} \right) f_{\mathcal{B}, \sigma_2} \left(  \lambda_{\mathcal{B}, i} \right)   \mathcal{P}_{\mathcal{B}, i} \nonumber \\
&=& F_{\mathcal{B}, \sigma_2} \star_N F_{\mathcal{B}, \sigma_1}. 
\end{eqnarray}

From the definition of $T_{\psi}$ provided by Eq.~\eqref{eq1:def pDTI} and $\psi_3 = \psi_1 \psi_2$, we have
\begin{eqnarray}
T_{\psi_3} &=& T_{\psi_1 \psi_2} = \int_{\Sigma_1 \times \Sigma_2} F_{\mathcal{A}, \sigma_1}  \star_N F_{\mathcal{A}, \sigma_2} \star_N \mathcal{X} \star_N   F_{\mathcal{B}, \sigma_2}  \star_N F_{\mathcal{B}, \sigma_1} d \left( \mu_1(\sigma_1) \times d \mu_2(\sigma_2) \right) \nonumber \\
&=& \int_{\Sigma_1} F_{\mathcal{A}, \sigma_1}  \star_N \left[  \int_{\Sigma_2}  F_{\mathcal{A}, \sigma_2} \star_N \mathcal{X} \star_N   F_{\mathcal{B}, \sigma_2}  d \mu_2(\sigma_2)  \right] \star_N F_{\mathcal{B}, \sigma_1} d \mu_1(\sigma_1)  \nonumber \\
&=& T_{\psi_1}\left( T_{\psi_2} (\mathcal{X}) \right).
\end{eqnarray}
Therefore, the mapping $\psi \rightarrow T_{\psi}$ is a hmomorphism.
$\hfill \Box$

\begin{lemma}\label{lma:power of 3-iii}
Let $f$ be a bounded real-valued function with the following properties for any given positive integer $m$:
\begin{eqnarray}
\phi_1 (\lambda_{\mathcal{A}}, \lambda_{\mathcal{B}}) = f(\lambda^{m}_{\mathcal{A}}), \mbox{and}~~\phi_2 (\lambda_{\mathcal{A}}, \lambda_{\mathcal{B}}) = f(\lambda^m_{\mathcal{B}}),
\end{eqnarray}
then 
\begin{eqnarray}
T_{\phi_1} (\mathcal{X}) = F^{m}_{\mathcal{A}} \star_N \mathcal{X}, \mbox{and}~~T_{\phi_2} (\mathcal{X}) =  \mathcal{X} \star_N F^{m}_{\mathcal{B}}, 
\end{eqnarray}
where 
\begin{eqnarray}
F^{m}_{\mathcal{A}} = \sum\limits_{i=1}^{\mathbb{I}_{1}^{N}} f(\lambda^{m}_{\mathcal{A},i}) \mathcal{P}_{\mathcal{A}, i}, \mbox{and}~~ F^{m}_{\mathcal{B}} = \sum\limits_{i=1}^{\mathbb{I}_{1}^{N}} f(\lambda^{m}_{\mathcal{B},i}) \mathcal{P}_{\mathcal{B}, i},
\end{eqnarray}
where $\mathcal{P}_{\mathcal{A}, i}$ and $\mathcal{P}_{\mathcal{B}, i}$ are the projection tensors for the underlying mappings $\phi_1 \rightarrow T_{\phi_1}$ and  $\phi_2 \rightarrow T_{\phi_2}$, respectively. We assume that $\sum\limits_{i=1}^{\mathbb{I}_{1}^{N}} f(\lambda^{m}_{\mathcal{A}}) \mathcal{P}_{\mathcal{A}, i}$ and $\sum\limits_{i=1}^{\mathbb{I}_{1}^{N}} f(\lambda^{m}_{\mathcal{B}}) \mathcal{P}_{\mathcal{B}, i}$ are positive definite tensors.  
\end{lemma}
\textbf{Proof:}
Since both functions $\phi_1$ and $\phi_2$ are belong to $\Phi$, this Lemma is proved by the definition of $T_{\psi}$ provided by Eq.~\eqref{eq1:def pDTI}  and Lemma~\ref{lma:kernel zero}.
$\hfill \Box$

We are ready to present the main theorem of this section. 

\begin{theorem}\label{thm:perturbation theory}
Let $f, g_{\mathcal{A}}, g_{\mathcal{B}}, h_{\mathcal{A}}, h_{\mathcal{B}}$ be bounded real-valued functions, and 
$\mathcal{E}_{\mathcal{A}}$ and $\mathcal{E}_{\mathcal{B}}$ be Hermitian tensors. We use $\mbox{Sp}(\mathcal{E}_{\mathcal{A}})$ and $\mbox{Sp}(\mathcal{E}_{\mathcal{B}})$ to represent the sets of eigenvalues of $\lambda_{\mathcal{A}}$ and $\lambda_{\mathcal{B}}$ for Hermitian tensors $\mathcal{E}_{\mathcal{A}}$ and $\mathcal{E}_{\mathcal{B}}$, respectively. We also assume that $m_{\mathcal{A}}, n_{\mathcal{A}},k_{\mathcal{A}}$ and $m_{\mathcal{B}}, n_{\mathcal{B}},k_{\mathcal{B}}$ are natural numbers. Let the function
\begin{eqnarray}\label{eq:1:thm:perturbation theory}
\psi(\lambda_{\mathcal{A}}, \lambda_{\mathcal{B}}) = \begin{cases}
 \frac{  h_{\mathcal{A}}(\lambda^{n_{\mathcal{A}}}_{\mathcal{A}})      }{  g_{\mathcal{A}}(\lambda^{m_{\mathcal{A}}}_{\mathcal{A}})   } \frac{    f \left( \lambda^{  k_{\mathcal{A}}  }_{\mathcal{A}} \right) -      f \left( \lambda^{  k_{\mathcal{B}}  }_{\mathcal{B}} \right)            }{ \lambda^{  k_{\mathcal{A}}  }_{\mathcal{A}}   - \lambda^{  k_{\mathcal{B}}  }_{\mathcal{B}}    } \frac{  h_{\mathcal{B}}(\lambda^{n_{\mathcal{B}}}_{\mathcal{B}})      }{  g_{\mathcal{B}}(\lambda^{m_{\mathcal{B}}}_{\mathcal{B}})   }, &  \mbox{if $(\lambda_{\mathcal{A}}, \lambda_{\mathcal{B}}) \in  \mbox{Sp}(\mathcal{E}_{\mathcal{A}}) \times \mbox{Sp}(\mathcal{E}_{\mathcal{B}})    $;}\\
0, & \mbox{otherwise.}
\end{cases}
\end{eqnarray}
Moreover, if we have
\begin{eqnarray}\label{eq:2:thm:perturbation theory}
G^{m_{\mathcal{A}}}_{\mathcal{A}} = g_{\mathcal{A}}(\mathcal{E}^{m_{\mathcal{A}}}_{\mathcal{A}}), ~~G^{m_{\mathcal{B}}}_{\mathcal{B}} = g_{\mathcal{B}}(\mathcal{E}^{m_{\mathcal{B}}}_{\mathcal{B}}) \nonumber \\
H^{n_{\mathcal{A}}}_{\mathcal{A}} = h_{\mathcal{A}}(\mathcal{E}^{n_{\mathcal{A}}}_{\mathcal{A}}), ~~H^{n_{\mathcal{B}}}_{\mathcal{B}} = h_{\mathcal{B}}(\mathcal{E}^{n_{\mathcal{B}}}_{\mathcal{B}}) \nonumber \\
F^{k_{\mathcal{A}}}_{\mathcal{A}} = f(\mathcal{E}^{k_{\mathcal{A}}}_{\mathcal{A}}), ~~F^{k_{\mathcal{B}}}_{\mathcal{B}} = f(\mathcal{E}^{k_{\mathcal{B}}}_{\mathcal{B}}),
\end{eqnarray}
then, 
\begin{eqnarray}\label{eq:3:thm:perturbation theory}
\mathcal{H}^{n_{\mathcal{A}}}_{\mathcal{A}} \star_N\left( F^{k_{\mathcal{A}}}_{\mathcal{A}} \star_N \mathcal{X} - \mathcal{X} \star_N F^{k_{\mathcal{B}}}_{\mathcal{B}}\right) \star_N \mathcal{H}^{n_{\mathcal{B}}}_{\mathcal{B}} = T_{\psi}\left( \mathcal{G}^{m_{\mathcal{A}}}_{\mathcal{A}} \star_N\left( E^{k_{\mathcal{A}}}_{\mathcal{A}} \star_N \mathcal{X} - \mathcal{X} \star_N E^{k_{\mathcal{B}}}_{\mathcal{B}}\right) \star_N \mathcal{G}^{m_{\mathcal{B}}}_{\mathcal{B}} \right).
\end{eqnarray}
In addition, we also have
\begin{eqnarray}\label{eq:4:thm:perturbation theory}
\lefteqn{\left\Vert \mathcal{H}^{n_{\mathcal{A}}}_{\mathcal{A}} \star_N\left( F^{k_{\mathcal{A}}}_{\mathcal{A}} \star_N \mathcal{X} - \mathcal{X} \star_N F^{k_{\mathcal{B}}}_{\mathcal{B}}\right) \star_N \mathcal{H}^{n_{\mathcal{B}}}_{\mathcal{B}} \right\Vert} \nonumber \\
&  &  \leq 
\left(\mathbb{I}_1^N \right)^2 \left\Vert \Psi \right\Vert   \left\Vert  \mathcal{G}^{m_{\mathcal{A}}}_{\mathcal{A}} \star_N\left( E^{k_{\mathcal{A}}}_{\mathcal{A}} \star_N \mathcal{X} - \mathcal{X} \star_N E^{k_{\mathcal{B}}}_{\mathcal{B}}\right) \star_N \mathcal{G}^{m_{\mathcal{B}}}_{\mathcal{B}}  \right\Vert 
\end{eqnarray}
\end{theorem}
\textbf{Proof:}

We define following functions with respect to $f, g_{\mathcal{A}}, g_{\mathcal{B}}, h_{\mathcal{A}}, h_{\mathcal{B}}$.
\begin{eqnarray}
\rho_{\mathcal{A}}(\lambda_{\mathcal{A}}, \lambda_{\mathcal{B}})  \define 
\lambda^{k_{\mathcal{A}}}_{\mathcal{A}} g_{\mathcal{A}}\left( \lambda^{m_{\mathcal{A}}}_{\mathcal{A}} \right)
g_{\mathcal{B}}\left( \lambda^{m_{\mathcal{B}}}_{\mathcal{B}} \right), ~~ \rho_{\mathcal{B}}(\lambda_{\mathcal{A}}, \lambda_{\mathcal{B}})  \define 
\lambda^{k_{\mathcal{B}}}_{\mathcal{B}} g_{\mathcal{A}}\left( \lambda^{m_{\mathcal{A}}}_{\mathcal{A}} \right)
g_{\mathcal{B}}\left( \lambda^{m_{\mathcal{B}}}_{\mathcal{B}} \right) \nonumber \\
\varsigma_{\mathcal{A}}(\lambda_{\mathcal{A}}, \lambda_{\mathcal{B}})  \define
f(\lambda^{k_{\mathcal{A}}}_{\mathcal{A}}) h_{\mathcal{A}}\left( \lambda^{n_{\mathcal{A}}}_{\mathcal{A}} \right)
h_{\mathcal{B}}\left( \lambda^{n_{\mathcal{B}}}_{\mathcal{B}} \right), ~~ \varsigma_{\mathcal{B}}(\lambda_{\mathcal{A}}, \lambda_{\mathcal{B}})  \define 
f(\lambda^{k_{\mathcal{B}}}_{\mathcal{B}}) h_{\mathcal{A}}\left( \lambda^{n_{\mathcal{A}}}_{\mathcal{A}} \right)
h_{\mathcal{B}}\left( \lambda^{n_{\mathcal{B}}}_{\mathcal{B}} \right),
\end{eqnarray}
where $(\lambda_{\mathcal{A}}, \lambda_{\mathcal{B}}) \in  \mbox{Sp}(\mathcal{E}_{\mathcal{A}}) \times \mbox{Sp}(\mathcal{E}_{\mathcal{B}})$. From Lemma~\ref{lma:power of 3-iii}, we have
\begin{eqnarray}
T_{\rho_{\mathcal{A}}}\left( \mathcal{X} \right) = G^{m_{\mathcal{A}}}_{\mathcal{A}} \star_N \mathcal{E}^{k_{\mathcal{A}}}_{\mathcal{A}} \star_N \mathcal{X} \star_N G^{m_{\mathcal{B}}}_{\mathcal{B}}, ~~
T_{\rho_{\mathcal{B}}}\left( \mathcal{X} \right) = G^{m_{\mathcal{A}}}_{\mathcal{A}} \star_N \mathcal{X} \star_N \mathcal{E}^{k_{\mathcal{B}}}_{\mathcal{B}} \star_N  G^{m_{\mathcal{B}}}_{\mathcal{B}}, \nonumber \\
T_{\varsigma_{\mathcal{A}}}\left( \mathcal{X} \right) = H^{n_{\mathcal{A}}}_{\mathcal{A}} \star_N \mathcal{F}^{k_{\mathcal{A}}}_{\mathcal{A}} \star_N \mathcal{X} \star_N H^{n_{\mathcal{B}}}_{\mathcal{B}}, ~~
T_{\varsigma_{\mathcal{B}}}\left( \mathcal{X} \right) = H^{n_{\mathcal{A}}}_{\mathcal{A}} \star_N \mathcal{X} \star_N \mathcal{F}^{k_{\mathcal{B}}}_{\mathcal{B}} \star_N  H^{n_{\mathcal{B}}}_{\mathcal{B}}.
\end{eqnarray}
By applying homomorphism of the mapping $\psi \rightarrow \Psi$ from Lemma~\ref{lma:homomorphism}, we have
\begin{eqnarray}
\lefteqn{T_{\psi}\left(  \mathcal{G}^{m_{\mathcal{A}}}_{\mathcal{A}} \star_N\left( E^{k_{\mathcal{A}}}_{\mathcal{A}} \star_N \mathcal{X} - \mathcal{X} \star_N E^{k_{\mathcal{B}}}_{\mathcal{B}}\right) \star_N \mathcal{G}^{m_{\mathcal{B}}}_{\mathcal{B}}  \right) } \nonumber \\
&& = T_{\psi}\left(T_{\rho_{\mathcal{A}}}\left(\mathcal{X}\right)  - 
T_{\rho_{\mathcal{B}}}\left(\mathcal{X}\right) \right)\nonumber \\
&& = T_{\psi\left( \rho_{\mathcal{A}} - \rho_{\mathcal{B}} \right)} \left(\mathcal{X}\right)= T_{\varsigma_{\mathcal{A}} - \varsigma_{\mathcal{B}}} \left(\mathcal{X} \right) \nonumber \\
&& = \mathcal{H}^{n_{\mathcal{A}}}_{\mathcal{A}} \star_N\left( F^{k_{\mathcal{A}}}_{\mathcal{A}} \star_N \mathcal{X} - \mathcal{X} \star_N F^{k_{\mathcal{B}}}_{\mathcal{B}}\right) \star_N \mathcal{H}^{n_{\mathcal{B}}}_{\mathcal{B}}. 
\end{eqnarray}
Therefore, Eq.~\eqref{eq:3:thm:perturbation theory} is established. 

Eq.~\eqref{eq:4:thm:perturbation theory} is true from Eq.~\eqref{eq:3:thm:perturbation theory} and Lemma~\ref{lma:T psi by psi}.
$\hfill \Box$

Following corollary is the variation of Theorem~\ref{thm:perturbation theory} by changing the negative sign in Eq.~\eqref{eq:1:thm:perturbation theory} to be the positive sign. The proof will be almost identical so we skip it.

\begin{corollary}\label{cor:perturbation theory}
Let $f, g_{\mathcal{A}}, g_{\mathcal{B}}, h_{\mathcal{A}}, h_{\mathcal{B}}$ be bounded real-valued functions, and 
$\mathcal{E}_{\mathcal{A}}$ and $\mathcal{E}_{\mathcal{B}}$ be Hermitian tensors. We use $\mbox{Sp}(\mathcal{E}_{\mathcal{A}})$ and $\mbox{Sp}(\mathcal{E}_{\mathcal{B}})$ to represent the sets of eigenvalues of $\lambda_{\mathcal{A}}$ and $\lambda_{\mathcal{B}}$ for Hermitian tensors $\mathcal{E}_{\mathcal{A}}$ and $\mathcal{E}_{\mathcal{B}}$, respectively. We also assume that $m_{\mathcal{A}}, n_{\mathcal{A}},k_{\mathcal{A}}$ and $m_{\mathcal{B}}, n_{\mathcal{B}},k_{\mathcal{B}}$ are natural numbers. Let the function
\begin{eqnarray}\label{eq:1:cor:perturbation theory}
\psi(\lambda_{\mathcal{A}}, \lambda_{\mathcal{B}}) = \begin{cases}
 \frac{  h_{\mathcal{A}}(\lambda^{n_{\mathcal{A}}}_{\mathcal{A}})      }{  g_{\mathcal{A}}(\lambda^{m_{\mathcal{A}}}_{\mathcal{A}})   } \frac{    f \left( \lambda^{  k_{\mathcal{A}}  }_{\mathcal{A}} \right)  +      f \left( \lambda^{  k_{\mathcal{B}}  }_{\mathcal{B}} \right)            }{ \lambda^{  k_{\mathcal{A}}  }_{\mathcal{A}}   + \lambda^{  k_{\mathcal{B}}  }_{\mathcal{B}}    } \frac{  h_{\mathcal{B}}(\lambda^{n_{\mathcal{B}}}_{\mathcal{B}})      }{  g_{\mathcal{B}}(\lambda^{m_{\mathcal{B}}}_{\mathcal{B}})   }, &  \mbox{if $(\lambda_{\mathcal{A}}, \lambda_{\mathcal{B}}) \in  \mbox{Sp}(\mathcal{E}_{\mathcal{A}}) \times \mbox{Sp}(\mathcal{E}_{\mathcal{B}})    $;}\\
0, & \mbox{otherwise.}
\end{cases}
\end{eqnarray}
Moreover, if we have
\begin{eqnarray}\label{eq:2:cor:perturbation theory}
G^{m_{\mathcal{A}}}_{\mathcal{A}} = g_{\mathcal{A}}(\mathcal{E}^{m_{\mathcal{A}}}_{\mathcal{A}}), ~~G^{m_{\mathcal{B}}}_{\mathcal{B}} = g_{\mathcal{B}}(\mathcal{E}^{m_{\mathcal{B}}}_{\mathcal{B}}) \nonumber \\
H^{n_{\mathcal{A}}}_{\mathcal{A}} = h_{\mathcal{A}}(\mathcal{E}^{n_{\mathcal{A}}}_{\mathcal{A}}), ~~H^{n_{\mathcal{B}}}_{\mathcal{B}} = h_{\mathcal{B}}(\mathcal{E}^{n_{\mathcal{B}}}_{\mathcal{B}}) \nonumber \\
F^{k_{\mathcal{A}}}_{\mathcal{A}} = f(\mathcal{E}^{k_{\mathcal{A}}}_{\mathcal{A}}), ~~F^{k_{\mathcal{B}}}_{\mathcal{B}} = f(\mathcal{E}^{k_{\mathcal{B}}}_{\mathcal{B}}),
\end{eqnarray}
then, 
\begin{eqnarray}\label{eq:3:cor:perturbation theory}
\mathcal{H}^{n_{\mathcal{A}}}_{\mathcal{A}} \star_N\left( F^{k_{\mathcal{A}}}_{\mathcal{A}} \star_N \mathcal{X} + \mathcal{X} \star_N F^{k_{\mathcal{B}}}_{\mathcal{B}}\right) \star_N \mathcal{H}^{n_{\mathcal{B}}}_{\mathcal{B}} = T_{\psi}\left( \mathcal{G}^{m_{\mathcal{A}}}_{\mathcal{A}} \star_N\left( E^{k_{\mathcal{A}}}_{\mathcal{A}} \star_N \mathcal{X} + \mathcal{X} \star_N E^{k_{\mathcal{B}}}_{\mathcal{B}}\right) \star_N \mathcal{G}^{m_{\mathcal{B}}}_{\mathcal{B}} \right).
\end{eqnarray}
In addition, we also have
\begin{eqnarray}\label{eq:4:cor:perturbation theory}
\lefteqn{\left\Vert \mathcal{H}^{n_{\mathcal{A}}}_{\mathcal{A}} \star_N\left( F^{k_{\mathcal{A}}}_{\mathcal{A}} \star_N \mathcal{X} + \mathcal{X} \star_N F^{k_{\mathcal{B}}}_{\mathcal{B}}\right) \star_N \mathcal{H}^{n_{\mathcal{B}}}_{\mathcal{B}} \right\Vert } \nonumber \\
&  & \leq 
\left( \mathbb{I}_1^N \right)^2 \left\Vert \Psi \right\Vert   \left\Vert  \mathcal{G}^{m_{\mathcal{A}}}_{\mathcal{A}} \star_N\left( E^{k_{\mathcal{A}}}_{\mathcal{A}} \star_N \mathcal{X} + \mathcal{X} \star_N E^{k_{\mathcal{B}}}_{\mathcal{B}}\right) \star_N \mathcal{G}^{m_{\mathcal{B}}}_{\mathcal{B}}  \right\Vert 
\end{eqnarray}
\end{corollary}

\section{Limiting Behavior of Random Parametrization Double Tensor Integrals}\label{sec:Limiting Behavior of Random Parametrizeation Double Tensor Integrals}

In this section, we will establish continuity of random PDTI. We need the following definition to define the convergence in mean for random tensors.

\begin{definition}\label{def:conv in mean}
We say that a sequence of random tensor $\mathcal{X}_n$ converges in the $r$-th mean towards the random tensor $\mathcal{X}$ with respect to the tensor norm $\left\Vert \cdot \right\Vert$, if we have
\begin{eqnarray}
\mathbb{E}\left( \left\Vert \mathcal{X}_n \right\Vert \right)~~~\mbox{exists,}
\end{eqnarray}
and
\begin{eqnarray}
\mathbb{E}\left( \left\Vert \mathcal{X} \right\Vert \right)~~~ \mbox{exists,}
\end{eqnarray}
and
\begin{eqnarray}
\lim\limits_{n \rightarrow \infty}\mathbb{E}\left( \left\Vert \mathcal{X}_n  - \mathcal{X} \right\Vert \right) = 0. 
\end{eqnarray}
We adopt the notatation $\mathcal{X}_n \xrightarrow[]{r} \mathcal{X}$ to represent that random tensors $\mathcal{X}_n$ converges in the $r$-th mean to the random tensor $\mathcal{X}$ with respect to the tensor norm $\left\Vert \cdot \right\Vert$.
\end{definition}

All limiting behaviors involving randomness discussed in this paper are based on convergence converges in the $1$-th mean.  

We define a special subset $\Psi_U$ within $\Psi$ that satisfies the following condition. If $\psi \in \Psi_U$, we have 
$\Sigma$, $f_{\mathcal{A}, \sigma}$ and $f_{\mathcal{B}, \sigma}$ in Eq.~\eqref{eq1:def pDTI} with the requirement that there is a increasing sequence of measurable subsets $S_k \in \Sigma$ for $i=1,2,\cdots$ such that 
\begin{eqnarray}\label{eq:Sigma decomp}
\Sigma = \bigcup\limits_{i=1}^{\infty}S_i,
\end{eqnarray}
and 
the family of functions $\{f_{\mathcal{A}, \sigma}, f_{\mathcal{B}, \sigma}\}$ is uniformly continous for every $i=1,2,\cdots$.

According to the $T_{\psi}(\mathcal{X})$ definition shown below, 
\begin{eqnarray}\label{eq1:def RPDTI}
T_{\psi}(\mathcal{X}) &=& \int_{\Sigma}  \left(\sum\limits_{i=1}^{\mathbb{I}_1^N}  f_{\mathcal{A}, \sigma} \left(  \lambda_{\mathcal{A}, i} \right)  \mathcal{P}_{\mathcal{A}, i} \right)\star_N \mathcal{X} \star_N \left(\sum\limits_{j=1}^{\mathbb{I}_1^N}  f_{\mathcal{B}, \sigma} \left(  \lambda_{\mathcal{B}, j} \right) \mathcal{P}_{\mathcal{B}, j} \right)  d \mu(\sigma),
\end{eqnarray}
the randomness of $T_{\psi}(\mathcal{X})$ comes from random variables $\lambda_{\mathcal{A}, i}, \lambda_{\mathcal{B}, i}$ and random tensors $\mathcal{P}_{\mathcal{A}, i}, \mathcal{P}_{\mathcal{B}, i}$.

\begin{lemma}\label{lma:conv of T psi t}
Let $\psi(\lambda_{\mathcal{A}}, \lambda_{\mathcal{B}}) = \int_{\Sigma} f_{\mathcal{A}, \sigma} (  \lambda_{\mathcal{A}} )  f_{\mathcal{B}, \sigma} (  \lambda_{\mathcal{B}} ) d \mu(\sigma) \in \Psi_U$ such that  functions $\left\Vert f^{(k)}_{\mathcal{A}, \sigma} \right\Vert_{\infty}$ and $\left\Vert f^{(k)}_{\mathcal{B}, \sigma} \right\Vert_{\infty}$ are bounded for $k=0,1,2$, where superscript $(k)$ represents the $k$-th derivative. The measure space $(\Sigma, \mu)$ follows Eq.~\ref{eq:Sigma decomp}. Also let $\{ \mathcal{E}_{\mathcal{A}, t} \}, \{ \mathcal{E}_{\mathcal{B}, t} \}$ be two indexed families of independent random Hermitian tensors for $t \in \mathbb{R}$ with formats~\footnote{In~\cite{potapov2010double}, their proof assumed that $\mathcal{E}_{\mathcal{A}, t} $ and $\mathcal{E}_{\mathcal{B}, t} $ is identical, but it can be more general by treating them as different objects.}: 
\begin{eqnarray}
\mathcal{E}_{\mathcal{A}, t} = \sum\limits_{i=1}^{\mathbb{I}_1^{N}} \lambda_{\mathcal{A}, i} \mathcal{P}_{\mathcal{A}, t, i},~~\mathcal{E}_{\mathcal{B}, t} = \sum\limits_{i=1}^{\mathbb{I}_1^{N}} \lambda_{\mathcal{B}, i} \mathcal{P}_{\mathcal{B}, t, i}
\end{eqnarray}
such that 
\begin{eqnarray}\label{eq1:lma:conv of T psi t}
\lim\limits_{t \rightarrow 0}\mathbb{E}\left( \left\Vert \mathcal{E}_{\mathcal{A}, t}   -  \mathcal{E}_{\mathcal{A}, 0}  \right\Vert\right) = 0, ~~ \lim\limits_{t \rightarrow 0}\mathbb{E}\left( \left\Vert \mathcal{E}_{\mathcal{B}, t}   -  \mathcal{E}_{\mathcal{B}, 0}  \right\Vert\right) = 0.
\end{eqnarray}
If $T_{\psi, t}$ is the random PDTI associated with $\psi$ and random tensors $ \mathcal{P}_{\mathcal{A}, t, i},  \mathcal{P}_{\mathcal{B}, t, i}$, then we have
\begin{eqnarray}
\lim\limits_{t \rightarrow 0}\mathbb{E}\left( \left\Vert T_{\psi, t} - T_{\psi, 0} \right\Vert\right) = 0. 
\end{eqnarray}
\end{lemma}
\textbf{Proof:}

Given $\epsilon > 0$, we wish to show that there is a function $\psi_{\epsilon} \in \Psi_U$ such that 
\begin{eqnarray}\label{eq2:lma:conv of T psi t}
\left\Vert \psi - \psi_{\epsilon} \right\Vert < \frac{\epsilon}{\left( \mathbb{I}_1^N \right)^2}.
\end{eqnarray}
The standard smoothing technique will be adopted here. We begin with the selecting the integer $i_{\epsilon} \in \mathbb{N}$ such that 
\begin{eqnarray}\label{eq3:lma:conv of T psi t}
\int_{\Sigma \backslash S_{i_{\epsilon}}} \left\Vert f_{\mathcal{A}, \sigma} (  \lambda_{\mathcal{A}} ) \right\Vert_{\infty}  \left\Vert f_{\mathcal{B}, \sigma} (  \lambda_{\mathcal{B}} ) \right\Vert_{\infty} d \mu(\sigma) < \frac{\epsilon}{3 \left( \mathbb{I}_1^N \right)^2}.
\end{eqnarray}
Then, given $ y_{\epsilon} > 0$, we set 
\begin{eqnarray}
f_{\mathcal{A}, \sigma, \epsilon} ( \lambda_{\mathcal{A}} )= \begin{cases}
f_{\mathcal{A}, \sigma} ( \lambda_{\mathcal{A}} ) \circledast \frac{  y_{\epsilon} }{\pi ( \lambda^2_{\mathcal{A}}+ y^2_{\epsilon})}, & \mbox{if $\sigma \in S_{i_{\epsilon}}$};\\
0, & \mbox{otherwise.}
\end{cases}
\end{eqnarray}
where $\circledast$ is the convolution operator. Similarly, we also set $f_{\mathcal{B}, \sigma, \epsilon} ( \lambda_{\mathcal{B}} )$ as 
\begin{eqnarray}
f_{\mathcal{B}, \sigma, \epsilon} ( \lambda_{\mathcal{B}} )= \begin{cases}
f_{\mathcal{B}, \sigma} ( \lambda_{\mathcal{B}} ) \circledast \frac{  y_{\epsilon} }{\pi ( \lambda^2_{\mathcal{B}}+ y^2_{\epsilon})}, & \mbox{if $\sigma \in S_{i_{\epsilon}}$};\\
0, & \mbox{otherwise.}
\end{cases}
\end{eqnarray}
By selecting the value $ y_{\epsilon}$ larger enought, we have 
\begin{eqnarray}\label{eq4:lma:conv of T psi t}
\left\Vert  f_{\mathcal{A}, \sigma, \epsilon} ( \lambda_{\mathcal{A}} ) -  f_{\mathcal{A}, \sigma} ( \lambda_{\mathcal{A}} ) \right\Vert_{\infty} <  \frac{\epsilon}{C_\epsilon \left( \mathbb{I}_1^N \right)^2},~~\mbox{and}~~ \left\Vert  f_{\mathcal{B}, \sigma, \epsilon} ( \lambda_{\mathcal{B}} ) -  f_{\mathcal{B}, \sigma} ( \lambda_{\mathcal{B}} )  \right\Vert_{\infty} < \frac{\epsilon}{C_\epsilon \left( \mathbb{I}_1^N \right)^2}.
\end{eqnarray}
where $C_{\epsilon}$ is defined as 
\begin{eqnarray}
C_{\epsilon} = 3 \mu(S_{i_\epsilon}) \max\limits_{\sigma \in S_{i_\epsilon}}\max(\left\Vert    f_{\mathcal{A}, \sigma} \right\Vert_{\infty},  \left\Vert    f_{\mathcal{B}, \sigma} \right\Vert_{\infty} ). 
\end{eqnarray}
The term $C_{\epsilon}$ is finite since the family of functions $\{f_{\mathcal{A}, \sigma}, f_{\mathcal{B}, \sigma}\}$ for $\sigma \in S_{i_\epsilon}$ is uniformly continous for every $i=1,2,\cdots$.
The function $\psi_{\epsilon}$ can be defined as 
\begin{eqnarray}
\psi_{\epsilon}(\lambda_{\mathcal{A}}, \lambda_{\mathcal{B}}) &\define& \int_{\Sigma} f_{\mathcal{A}, \sigma, \epsilon} (  \lambda_{\mathcal{A}} )  f_{\mathcal{B}, \sigma, \epsilon} (  \lambda_{\mathcal{B}} ) d \mu(\sigma)
\nonumber \\
&=& \int_{S_{i_{\epsilon}}} f_{\mathcal{A}, \sigma, \epsilon} (  \lambda_{\mathcal{A}} )  f_{\mathcal{B}, \sigma, \epsilon} (  \lambda_{\mathcal{B}} ) d \mu(\sigma),
\end{eqnarray}
then, we have
\begin{eqnarray}\label{eq5:lma:conv of T psi t}
\psi(\lambda_{\mathcal{A}}, \lambda_{\mathcal{B}}) - \psi_{\epsilon}(\lambda_{\mathcal{A}}, \lambda_{\mathcal{B}}) &=& \int_{S_{i_{\epsilon}}} f_{\mathcal{A}, \sigma} (  \lambda_{\mathcal{A}} )  \left[ f_{\mathcal{B}, \sigma} (  \lambda_{\mathcal{B}} ) -   f_{\mathcal{B}, \sigma, \epsilon} (  \lambda_{\mathcal{B}} )  \right] d \mu(\sigma)  \nonumber \\
& & +  \int_{S_{i_{\epsilon}}}   \left[ f_{\mathcal{A}, \sigma} (  \lambda_{\mathcal{A}} ) -   f_{\mathcal{A}, \sigma, \epsilon} (  \lambda_{\mathcal{A}} )  \right]  f_{\mathcal{B}, \sigma, \epsilon} (  \lambda_{\mathcal{B}} ) d \mu(\sigma) \nonumber \\
& & +  \int_{\Sigma \backslash S_{i_{\epsilon}}}    f_{\mathcal{A}, \sigma} (  \lambda_{\mathcal{A}} )     f_{\mathcal{B}, \sigma, \epsilon} (  \lambda_{\mathcal{B}} ) d \mu(\sigma)
\end{eqnarray}
By applying Eqs.~\eqref{eq3:lma:conv of T psi t} and~\eqref{eq4:lma:conv of T psi t} to Eq.~\eqref{eq5:lma:conv of T psi t}, we can have $\left\Vert \psi - \psi_{\epsilon} \right\Vert_{\Psi} <  \frac{\epsilon}{\left( \mathbb{I}_1^N \right)^2}$. 

Our next goal is to show 
\begin{eqnarray}\label{eq6:lma:conv of T psi t}
\left\Vert T_{\psi_{\epsilon}, t}(\mathcal{X}) - T_{\psi_{\epsilon}, 0}(\mathcal{X}) \right\Vert <  \left\Vert \mathcal{X} \right\Vert \epsilon ,~~\mbox{if $t  < \delta$.} 
\end{eqnarray}

If we set
\begin{eqnarray}
F_{\mathcal{A}, \sigma, t, \epsilon} \define f_{\mathcal{A}, \sigma, \epsilon}\left( \mathcal{E}_{\mathcal{A}, t} \right) = \sum\limits_{i=1}^{\mathbb{I}_1^{N}} f_{\mathcal{A}, \sigma, \epsilon}\left( \lambda_{\mathcal{A}, i} \right) \mathcal{P}_{\mathcal{A}, t, i},
\end{eqnarray}
and 
\begin{eqnarray}
F_{\mathcal{B}, \sigma, t, \epsilon} \define f_{\mathcal{B}, \sigma, \epsilon}\left( \mathcal{E}_{\mathcal{B}, t} \right) = \sum\limits_{i=1}^{\mathbb{I}_1^{N}} f_{\mathcal{B}, \sigma, \epsilon}\left( \lambda_{\mathcal{B}, i} \right) \mathcal{P}_{\mathcal{B}, t, i};
\end{eqnarray}
then, from Theorem~\ref{thm:perturbation theory} and Theorem 4  of~\cite{potapov2009unbounded},  we have 
\begin{eqnarray}\label{eq6-1:lma:conv of T psi t}
\left\Vert   F_{\mathcal{A}, \sigma, t, \epsilon} - F_{\mathcal{A}, \sigma, 0, \epsilon}    \right\Vert \leq 
\overline{C}_{\epsilon} \left\Vert \mathcal{E}_{\mathcal{A}, t} - \mathcal{E}_{\mathcal{A}, 0}       \right\Vert,
\end{eqnarray}
and
\begin{eqnarray}\label{eq6-2:lma:conv of T psi t}
\left\Vert   F_{\mathcal{B}, \sigma, t, \epsilon} - F_{\mathcal{B}, \sigma, 0, \epsilon}    \right\Vert \leq 
\overline{C}_{\epsilon} \left\Vert \mathcal{E}_{\mathcal{B}, t} - \mathcal{E}_{\mathcal{B}, 0}       \right\Vert,
\end{eqnarray}
where the contant $\overline{C}_{\epsilon}$ can be expressed as
\begin{eqnarray}
\overline{C}_{\epsilon} = \max\limits_{\sigma \in S_{i_\epsilon}}\max \left(   \max\limits_{i=0,1,2}\left\Vert    f^{(i)}_{\mathcal{A}, \sigma} \right\Vert_{\infty},  \max\limits_{i=0,1,2} \left\Vert    f^{(i)}_{\mathcal{B}, \sigma} \right\Vert_{\infty} \right),
\end{eqnarray}
where the superscript $(i)$ is the $i$-th derivative.

By taking expectations for the both sides of Eqs.~\eqref{eq6-1:lma:conv of T psi t} and~\eqref{eq6-2:lma:conv of T psi t}, and from the assumptions provided by Eq.~\eqref{eq1:lma:conv of T psi t}, we have
\begin{eqnarray}\label{eq7:lma:conv of T psi t}
\mathbb{E} \left(\left\Vert   F_{\mathcal{A},\sigma, t, \epsilon} - F_{\mathcal{A},\sigma, 0, \epsilon}    \right\Vert \right) \leq \left(\frac{\epsilon}{\mu (S_{i_\epsilon})}\right)^{1/2}~~\mbox{, if $t  < \delta$ such that $ \mathbb{E}\left(\left\Vert \mathcal{E}_{\mathcal{A}, t} - \mathcal{E}_{\mathcal{A}, 0} \right\Vert \right)   < \frac{ \sqrt{\epsilon} }{\sqrt{\mu(S_{i_\epsilon})} \overline{C}_{\epsilon} }$,} 
\end{eqnarray}
where $\mu (S_{i_\epsilon})$ is the measure for the $S_{i_\epsilon}$. Similarly, we also have
\begin{eqnarray}\label{eq8:lma:conv of T psi t}
\mathbb{E} \left(\left\Vert   F_{\mathcal{B},\sigma, t, \epsilon} - F_{\mathcal{B},\sigma, 0, \epsilon}    \right\Vert \right) \leq \left(\frac{\epsilon}{\mu (S_{i_\epsilon})}\right)^{1/2}~~\mbox{, if $t  < \delta$ such that $ \mathbb{E}\left(\left\Vert \mathcal{E}_{\mathcal{A}, t} - \mathcal{E}_{\mathcal{A}, 0} \right\Vert \right)   < \frac{ \sqrt{\epsilon}  }{\sqrt{\mu(S_{i_\epsilon})} \overline{C}_{\epsilon} }$.} 
\end{eqnarray}
Then, we have
\begin{eqnarray}\label{eq9:lma:conv of T psi t}
\mathbb{E}\left(\left\Vert T_{\psi_{\epsilon}, t}(\mathcal{X}) - T_{\psi_{\epsilon}, 0}(\mathcal{X}) \right\Vert \right)
&=& \mathbb{E}\left( \left\Vert \int_{S_{i_\epsilon}} \left[ F_{\mathcal{A},\sigma, t, \epsilon} - F_{\mathcal{A},\sigma, 0, \epsilon}  \right] \star_N \mathcal{X} \star_N   \left[ F_{\mathcal{B},\sigma, t, \epsilon} - F_{\mathcal{B},\sigma, 0, \epsilon}  \right] d \mu(\sigma) \right\Vert \right) \nonumber \\
&\leq & \left[ \int_{S_{i_\epsilon}}  \mathbb{E}\left( F_{\mathcal{A},\sigma, t, \epsilon} - F_{\mathcal{A},\sigma, 0, \epsilon}    \right)    \mathbb{E}\left( F_{\mathcal{B},\sigma, t, \epsilon} - F_{\mathcal{B},\sigma, 0, \epsilon}    \right) d \mu(\sigma) \right] \left\Vert \mathcal{X} \right\Vert
\end{eqnarray}
By taking expectation of the both sides of Eq.~\eqref{eq9:lma:conv of T psi t} and applying Eqs.~\eqref{eq7:lma:conv of T psi t} and~\eqref{eq8:lma:conv of T psi t}, we obtain
\begin{eqnarray}\label{eq10:lma:conv of T psi t}
\mathbb{E}\left( \left\Vert T_{\psi_{\epsilon}, t}(\mathcal{X}) - T_{\psi_{\epsilon}, 0}(\mathcal{X}) \right\Vert \right)
\leq  \left\Vert \mathcal{X} \right\Vert  \epsilon .
\end{eqnarray}

Finally, given $t < \delta$, we have
\begin{eqnarray}
\mathbb{E}\left(\left\Vert T_{\psi, t}(\mathcal{X}) - T_{\psi, 0}(\mathcal{X}) \right\Vert \right)
&\leq & \mathbb{E}\left(\left\Vert T_{\psi, t}(\mathcal{X}) - T_{\psi_{\epsilon}, t}(\mathcal{X}) \right\Vert \right)
+ \mathbb{E}\left(\left\Vert T_{\psi_{\epsilon}, t}(\mathcal{X}) - T_{\psi_{\epsilon}, 0}(\mathcal{X}) \right\Vert \right) \nonumber \\
&  &+ \mathbb{E}\left(\left\Vert T_{\psi_{\epsilon}, 0}(\mathcal{X}) - T_{\psi, 0}(\mathcal{X}) \right\Vert \right) \nonumber \\
&\leq& 3  \left\Vert \mathcal{X} \right\Vert \epsilon,
\end{eqnarray} 
where the first and third terms are obtained from Eq.~\eqref{eq2:lma:conv of T psi t} and Lemma~\ref{lma:T psi by psi}, and the second term comes from Eq~\eqref{eq10:lma:conv of T psi t}. 
$\hfill \Box$

Following Lemma is the derivative tensor relation after the action of $T_{\psi}$.

\begin{lemma}\label{lma:derivative after T psi action}
Let $\mathcal{E}_{t}$ for $t \in \mathbb{R}$ be a family of Hermitian tensors such that 
\begin{eqnarray}\label{eq1:lma:derivative after T psi action}
\lim\limits_{t \rightarrow 0} \left\Vert \mathcal{E}_{t}   -  \mathcal{E}_{0}  \right\Vert= 0.
\end{eqnarray} 
Moreover, we also have
\begin{eqnarray}\label{eq2:lma:derivative after T psi action}
G^{m}_{0} = g(\mathcal{E}^{m}_{0}), ~~G^{m}_{t} = g(\mathcal{E}^{m}_{t}); \nonumber \\
H^{n}_{0} = h(\mathcal{E}^{n}_{0}), ~~H^{n}_{t} = h(\mathcal{E}^{n}_{t});
\nonumber \\
F^{k}_{0} = f(\mathcal{E}^{k}_{0}), ~~F^{k}_{t} = f(\mathcal{E}^{k}_{t}).
\end{eqnarray}

If $\psi \in \Psi_U$ and if 
\begin{eqnarray} 
\mathcal{A}_0 = \lim\limits_{t \rightarrow 0} G^{m}_{t}\frac{\mathcal{E}^k_{t} - \mathcal{E}^k_{0}}{t^k}  G^{m}_{0}
\end{eqnarray}
exists, then the limit 
\begin{eqnarray} 
\mathcal{B}_0 = \lim\limits_{t \rightarrow 0} H^{n}_{t}\frac{\mathcal{F}^k_{t} - \mathcal{F}^k_{0}}{t^k}  H^{n}_{0}
\end{eqnarray}
exist. Moreover, we have 
\begin{eqnarray} 
\mathcal{B}_0 = T_{\psi}\left( \mathcal{A}_0 \right),
\end{eqnarray}
where $\psi$ can be expressed as 
\begin{eqnarray}\label{eq3:lma:derivative after T psi action}
\psi(\lambda_{\mathcal{A}}, \lambda_{\mathcal{B}}) = \begin{cases}
 \frac{  h(\lambda^{n}_{\mathcal{A}})      }{  g(\lambda^{m}_{\mathcal{A}})   } \frac{    f \left( \lambda^{  k  }_{\mathcal{A}} \right) -      f \left( \lambda^{  k }_{\mathcal{B}} \right)            }{ \lambda^{  k }_{\mathcal{A}}   - \lambda^{  k }_{\mathcal{B}}    } \frac{  h(\lambda^{n}_{\mathcal{B}})      }{  g(\lambda^{m}_{\mathcal{B}})   }, &  \mbox{if $(\lambda_{\mathcal{A}}, \lambda_{\mathcal{B}}) \in  \mbox{Sp}(\mathcal{E}_{0}) \times \mbox{Sp}(\mathcal{E}_{0})    $;}\\
0, & \mbox{otherwise.}
\end{cases}
\end{eqnarray}
\end{lemma}
\textbf{Proof:}

By setting 
\begin{eqnarray} 
\mathcal{A}_t = G^{m}_{t}\frac{\mathcal{E}^k_{t} - \mathcal{E}^k_{0}}{t^k}  G^{m}_{0},
\end{eqnarray}
and
\begin{eqnarray} 
\mathcal{B}_t= H^{n}_{t}\frac{\mathcal{F}^k_{t} - \mathcal{F}^k_{0}}{t^k}  H^{n}_{0},
\end{eqnarray}
we have $\mathcal{B}_t = T_{\psi, t}(\mathcal{A}_t)$ from Theorem~\ref{thm:perturbation theory}. Then, we have
\begin{eqnarray} 
\lim\limits_{t \rightarrow 0}\left\Vert \mathcal{B}_t - \mathcal{B}_0 \right\Vert &\leq&
 \lim\limits_{t \rightarrow 0}\left( \left\Vert T_{\psi, t}(\mathcal{A}_t - \mathcal{A}_0 )\right\Vert +
\left\Vert (T_{\psi, t} - T_{\psi}) (\mathcal{A}_0)\right\Vert  \right) \nonumber \\
& \leq  & \epsilon,
\end{eqnarray}
where $\lim\limits_{t \rightarrow 0} \mathcal{A}_t= \mathcal{A}_0$ comes from the assumption provided by Eq.~\eqref{eq1:lma:derivative after T psi action}, and $T_{\psi, t} - T_{\psi}$ comes from Lemma~\ref{lma:conv of T psi t}.
$\hfill \Box$

If we have the following condition in Lemma~\ref{lma:derivative after T psi action}, $\mathcal{E}_{t}$ for $t \in \mathbb{R}$ be a family of random  Hermitian tensors such that 
\begin{eqnarray}\label{eq1:lma:derivative after T psi action r}
\lim\limits_{t \rightarrow 0}\mathbb{E}\left( \left\Vert \mathcal{E}_{t}   -  \mathcal{E}_{0}  \right\Vert \right)= 0,
\end{eqnarray} 
then, we have
\begin{eqnarray} 
\lim\limits_{t \rightarrow 0}\mathbb{E}\left( \left\Vert T_{\psi, t}\left( \mathcal{A}_t \right) - \mathcal{B}_0 \right\Vert  \right) = 0.
\end{eqnarray}
The proof is similar to Lemma~\ref{lma:derivative after T psi action}.

\section{New Inequalities By PDTI}\label{sec:New Inequalities By PDTI} 

In this section, we will apply the proposed PDTI to derive several new inequalities.

\begin{theorem}\label{thm:11lma}
Let $\mathcal{A}, \mathcal{B} \in \mathbb{C}^{I_1 \times \cdots \times I_N \times I_1 \times \cdots \times I_N} $ be random Hermitian tensors and $\mathcal{X} \in \mathbb{C}^{I_1 \times \cdots \times I_N \times I_1 \times \cdots \times I_N} $ be a Hermitian tensor. For every $0 \leq \omega \leq m$, we have
\begin{eqnarray}\label{eq1:thm:11lma}
\mathrm{Pr}\left( \left\Vert \mathcal{A}^m \star_N \mathcal{X} \star_N   \mathcal{B}^{\omega} -  \mathcal{A}^{\omega} \star_N \mathcal{X} \star_N   \mathcal{B}^{m}  \right\Vert \geq \theta \right) &\leq&  
\frac{    \left( \mathbb{I}_1^N \right)^2 }{\theta}  \left[\min\limits_{c>0} \left( \sqrt{2c}\left\Vert g(t) \right\Vert_2 + \sqrt{2/c}\left\Vert g'(t) \right\Vert_2    \right) \right] \nonumber \\
 &  & \times   \mathbb{E}\left( \left\Vert \mathcal{A}^m \star_N \mathcal{X} - \mathcal{X} \star_N \mathcal{B}^m \right\Vert \right),
 \end{eqnarray}
where $g(t)$ is
\begin{eqnarray}\label{eq2:thm:11lma}
g(t) \define  \frac{\exp\left(\frac{(m- 2 \omega ) t}{2}\right) -  \exp\left(\frac{(2 \omega - m ) t}{2}\right)    }{\exp\left( \frac{mt}{2} \right) -  \exp\left(\frac{-mt}{2} \right)   }
\end{eqnarray}
\end{theorem}
\textbf{Proof:}

From Theorem~\ref{thm:perturbation theory}, we have 
\begin{eqnarray}\label{eq3:thm:11lma}
\mathcal{A}^m \star_N \mathcal{X} \star_N   \mathcal{B}^{\omega} -  \mathcal{A}^{\omega} \star_N \mathcal{X} \star_N   \mathcal{B}^{m} =T_{\psi}\left(  \mathcal{A}^m \star_N \mathcal{X} - \mathcal{X} \star_N \mathcal{B}^m\right),
\end{eqnarray} 
where $\psi$ is 
\begin{eqnarray}\label{eq4:thm:11lma}
\psi(\lambda_{\mathcal{A}}, \lambda_{\mathcal{B}}) =  \frac{ \lambda^{m - \omega}_{\mathcal{A}}\lambda^{  \omega}_{\mathcal{B}} - \lambda^{ \omega}_{\mathcal{A}}\lambda^{ m - \omega}_{\mathcal{B}}           }{ \lambda^m_{\mathcal{A}}  - \lambda^m_{\mathcal{B}} }.
\end{eqnarray}
Then, Eq.~\eqref{eq4:thm:11lma} will be obtained by setting $t = \log \frac{\lambda_{\mathcal{A}}}{ \lambda_{\mathcal{B}}}$ in Eq.~\eqref{eq2:thm:11lma}.

By applying Lemma~\ref{lma:T psi by psi} and Corollary~\ref{cor:bound for psi norm Fourier} to the function $g(t)$ provided by Eq.~\eqref{eq2:thm:11lma}, we have
\begin{eqnarray}
 \left\Vert \mathcal{A}^m \star_N \mathcal{X} \star_N   \mathcal{B}^{\omega} -  \mathcal{A}^{\omega} \star_N \mathcal{X} \star_N   \mathcal{B}^{m}  \right\Vert 
&\leq&  \left( \mathbb{I}_1^N \right)^2  \left[\min\limits_{c>0} \left( \sqrt{2c}\left\Vert g(t) \right\Vert_2 + \sqrt{2/c}\left\Vert g'(t) \right\Vert_2    \right) \right] \nonumber \\
 &  & \times   \left\Vert \mathcal{A}^m \star_N \mathcal{X} - \mathcal{X} \star_N \mathcal{B}^m \right\Vert.
\end{eqnarray}
Therefore, we have
\begin{eqnarray}\label{eq5:thm:11lma}
\lefteqn{\mathrm{Pr}\left( \left\Vert \mathcal{A}^m \star_N \mathcal{X} \star_N   \mathcal{B}^{\omega} -  \mathcal{A}^{\omega} \star_N \mathcal{X} \star_N   \mathcal{B}^{m}  \right\Vert \geq \theta \right) } \nonumber \\
&\leq&  
\mathrm{Pr}\left( \left\{ \left( \mathbb{I}_1^N \right)^2  \left[\min\limits_{c>0} \left( \sqrt{2c}\left\Vert g(t) \right\Vert_2 + \sqrt{2/c}\left\Vert g'(t) \right\Vert_2    \right) \right]   \left\Vert \mathcal{A}^m \star_N \mathcal{X} - \mathcal{X} \star_N \mathcal{B}^m \right\Vert \right\} \geq \theta \right) \nonumber \\
&=& \mathrm{Pr}\left(    \left\Vert \mathcal{A}^m \star_N \mathcal{X} - \mathcal{X} \star_N \mathcal{B}^m \right\Vert  \geq \frac{\theta}{    \left( \mathbb{I}_1^N \right)^2  \left[\min\limits_{c>0} \left( \sqrt{2c}\left\Vert g(t) \right\Vert_2 + \sqrt{2/c}\left\Vert g'(t) \right\Vert_2    \right) \right]       }     \right).
\end{eqnarray}
This theorem is proved by applying Markov inequality to Eq.~\eqref{eq5:thm:11lma}.
$\hfill \Box$

If $m=1$, Theorem~\ref{thm:11lma} becomes the tail bound for Heinz inequality~\cite{kosaki2011positive}.

Following corollary is obtained by applying Corollary~\ref{cor:perturbation theory} to the same conditions of Theorem~\ref{thm:11lma} for the tensor $ \mathcal{A}^m \star_N \mathcal{X} \star_N   \mathcal{B}^{\omega} + \mathcal{A}^{\omega} \star_N \mathcal{X} \star_N   \mathcal{B}^{m}$. We will skip the proof here due to the similarity of the proof provided by Theorem~\ref{thm:11lma}.

\begin{corollary}\label{cor:11lma}
Let $\mathcal{A}, \mathcal{B} \in \mathbb{C}^{I_1 \times \cdots \times I_N \times I_1 \times \cdots \times I_N} $ be random Hermitian tensors and $\mathcal{X} \in \mathbb{C}^{I_1 \times \cdots \times I_N \times I_1 \times \cdots \times I_N} $ be a Hermitian tensor. For every $0 \leq \omega \leq m$, we have
\begin{eqnarray}\label{eq1:cor:11lma}
\mathrm{Pr}\left( \left\Vert \mathcal{A}^m \star_N \mathcal{X} \star_N   \mathcal{B}^{\omega} + \mathcal{A}^{\omega} \star_N \mathcal{X} \star_N   \mathcal{B}^{m}  \right\Vert \geq \theta \right) &\leq&  
\frac{    \left( \mathbb{I}_1^N \right)^2 }{\theta}  \left[\min\limits_{c>0} \left( \sqrt{2c}\left\Vert g(t) \right\Vert_2 + \sqrt{2/c}\left\Vert g'(t) \right\Vert_2    \right) \right] \nonumber \\
 &  & \times   \mathbb{E}\left( \left\Vert \mathcal{A}^m \star_N \mathcal{X} + \mathcal{X} \star_N \mathcal{B}^m \right\Vert \right)
 \end{eqnarray}
where $g(t)$ is
\begin{eqnarray}\label{eq2:cor:11lma}
g(t) \define  \frac{\exp\left(\frac{(m- 2 \theta) t}{2}\right) +  \exp\left(\frac{(2 \theta - m ) t}{2}\right)    }{\exp\left( \frac{mt}{2} \right) +  \exp\left(\frac{-mt}{2} \right)   }
\end{eqnarray}
\end{corollary}

Before presenting the following theorem, we have to introduce some notations. Given the tensor $\mathcal{A} \in \mathbb{C}^{I_1 \times \cdots \times I_N \times I_1 \times \cdots \times I_N}$, we use the absolute symbol $\left\vert \mathcal{A}\right\vert$ to represent the following:
\begin{eqnarray}
\left\vert \mathcal{A}\right\vert \define \sqrt{\mathcal{A}^{\mathrm{T}} \star_N \mathcal{A}}.
\end{eqnarray} 
Also, we use the symbol $[\mathcal{A}, \mathcal{B}]$, where $\mathcal{A}, \mathcal{B} \in \mathbb{C}^{I_1 \times \cdots \times I_N \times I_1 \times \cdots \times I_N}$, to represent the commutator between two tensors, it is defined as:
\begin{eqnarray}
[\mathcal{A}, \mathcal{B}] \define \mathcal{A} \star_N \mathcal{B} - \mathcal{B} \star_N \mathcal{A}.
\end{eqnarray}

\begin{theorem}\label{thm:15lma}
Let $\mathcal{A}, \mathcal{B} \in \mathbb{C}^{I_1 \times \cdots \times I_N \times I_1 \times \cdots \times I_N} $ be random positive definite tensors and $\mathcal{X} \in \mathbb{C}^{I_1 \times \cdots \times I_N \times I_1 \times \cdots \times I_N} $ be a Hermitian tensor. For every $0 \leq \nu \leq 1$ and two nonnegative real numbers $r_0, r_1$ satisfying $r_0 + r_1 = 1$, we have
\begin{eqnarray}\label{eq1:thm:15lma}
\mathrm{Pr}\left( \left\Vert [\mathcal{A}  \left\vert \mathcal{A} \right\vert^{-\nu}, \mathcal{B}]  \right\Vert \geq \theta \right) &\leq&  
\frac{    \left( \mathbb{I}_1^N \right)^2 }{\theta}  \left[\min\limits_{c>0} \left( \sqrt{2c}\left\Vert g(t) \right\Vert_2 + \sqrt{2/c}\left\Vert g'(t) \right\Vert_2    \right) \right] \nonumber \\
 &  & \times   \mathbb{E}\left( \left\Vert \left\vert \mathcal{A} \right\vert^{- r_0 \nu } \star_N [\mathcal{A}, \mathcal{B}] \star_N  \left\vert \mathcal{A} \right\vert^{- r_1 \nu}    \right\Vert \right)
 \end{eqnarray}
where $g(t)$ is
\begin{eqnarray}\label{eq2:thm:15lma}
g(t) \define  \frac{\exp\left[\frac{(1 - 2 r_1 \nu) t}{2}\right] -  \exp\left[\frac{(2r_0 \nu -1 ) t}{2}\right]    }{\exp\left( \frac{t}{2} \right) -  \exp\left(\frac{-t}{2} \right)   }.
\end{eqnarray}
\end{theorem}
\textbf{Proof:}

From Theorem~\ref{thm:perturbation theory}, we have 
\begin{eqnarray}\label{eq3:thm:15lma}
\left[\mathcal{A}  \left\vert \mathcal{A} \right\vert^{-\nu}, \mathcal{B} \right] =T_{\psi}\left(  \left\vert \mathcal{A} \right\vert^{-r_0 \nu} \star_N [\mathcal{A}, \mathcal{B}] \star_N \left\vert \mathcal{A} \right\vert^{-r_1\nu }  \right),
\end{eqnarray} 
where $\psi$ is 
\begin{eqnarray}\label{eq4:thm:15lma}
\psi(\lambda_{\mathcal{A}}, \lambda_{\mathcal{B}}) = \lambda^{r_0 \nu}_{\mathcal{A}} \frac{  \lambda^{1 - \nu}_{\mathcal{A}}  - \lambda^{1 - \nu}_{\mathcal{B}}        }{ \lambda_{\mathcal{A}}  - \lambda_{\mathcal{B}} }\lambda^{r_1 \nu }_{\mathcal{B}}  .
\end{eqnarray}
Then, Eq.~\eqref{eq2:thm:15lma} will be obtained by setting $t = \log \frac{\lambda_{\mathcal{A}}}{ \lambda_{\mathcal{B}}}$ in Eq.~\eqref{eq4:thm:15lma}.

By applying Lemma~\ref{lma:T psi by psi} and Corollary~\ref{cor:bound for psi norm Fourier} to the function $g(t)$ provided by Eq.~\eqref{eq2:thm:15lma}, we have
\begin{eqnarray}
 \left\Vert  \left[\mathcal{A}  \left\vert \mathcal{A} \right\vert^{-\nu}, \mathcal{B} \right] \right\Vert 
&\leq&  \left( \mathbb{I}_1^N \right)^2  \left[\min\limits_{c>0} \left( \sqrt{2c}\left\Vert g(t) \right\Vert_2 + \sqrt{2/c}\left\Vert g'(t) \right\Vert_2    \right) \right] \nonumber \\
 &  & \times   \left\Vert \left\vert \mathcal{A} \right\vert^{-r_0 \nu} \star_N [\mathcal{A}, \mathcal{B}] \star_N \left\vert \mathcal{A} \right\vert^{-r_1\nu }  \right\Vert.
\end{eqnarray}
Therefore, we have
\begin{eqnarray}\label{eq5:thm:15lma}
\lefteqn{\mathrm{Pr}\left( \left\Vert \left[\mathcal{A}  \left\vert \mathcal{A} \right\vert^{-\nu}, \mathcal{B} \right]  \right\Vert \geq \theta \right) } \nonumber \\
&\leq&  
\mathrm{Pr}\left( \left\{ \left( \mathbb{I}_1^N \right)^2  \left[\min\limits_{c>0} \left( \sqrt{2c}\left\Vert g(t) \right\Vert_2 + \sqrt{2/c}\left\Vert g'(t) \right\Vert_2    \right) \right] \left\Vert \left\vert \mathcal{A} \right\vert^{-r_0 \nu} \star_N [\mathcal{A}, \mathcal{B}] \star_N \left\vert \mathcal{A} \right\vert^{-r_1\nu }  \right\Vert \right\} \geq \theta \right) \nonumber \\
&=& \mathrm{Pr}\left(    \left\Vert \left\vert \mathcal{A} \right\vert^{-r_0 \nu} \star_N [\mathcal{A}, \mathcal{B}] \star_N \left\vert \mathcal{A} \right\vert^{-r_1\nu }  \right\Vert \geq \frac{\theta}{    \left( \mathbb{I}_1^N \right)^2  \left[\min\limits_{c>0} \left( \sqrt{2c}\left\Vert g(t) \right\Vert_2 + \sqrt{2/c}\left\Vert g'(t) \right\Vert_2    \right) \right]       }     \right).
\end{eqnarray}
This theorem is proved by applying Markov inequality to Eq.~\eqref{eq5:thm:15lma}.
$\hfill \Box$

\begin{theorem}\label{thm:16thm}
Let $\mathcal{A}, \mathcal{B} \in \mathbb{C}^{I_1 \times \cdots \times I_N \times I_1 \times \cdots \times I_N} $ be random positive definite tensors. For every $0 \leq \omega \leq 1$ and positive integers $m,n$, we have
\begin{eqnarray}\label{eq1:thm:16thm}
\mathrm{Pr}\left( \left\Vert \mathcal{A}^{n\omega} - \mathcal{B}^{m\omega} \right\Vert \geq \theta \right) &\leq&  
\frac{    \left( \mathbb{I}_1^N \right)^2 }{\omega \theta}  \left[\min\limits_{c>0} \left( \sqrt{2c}\left\Vert g(t) \right\Vert_2 + \sqrt{2/c}\left\Vert g'(t) \right\Vert_2    \right) \right] \nonumber \\
 &  & \times   \mathbb{E}\left( \left\Vert  \mathcal{A}^{n} - \mathcal{B}^{m}  \right\Vert^{\omega} \right)
 \end{eqnarray}
where $g(t)$ is
\begin{eqnarray}\label{eq2:thm:16thm}
g(t) \define  \frac{\exp\left( \frac{\omega t}{2}\right) -  \exp\left(\frac{-\omega t}{2}\right)    }{\exp\left( \frac{t}{2} \right) -  \exp\left(\frac{-t}{2} \right)   }.
\end{eqnarray}
\end{theorem}
\textbf{Proof:}
Since the spectral norm is the same for taking a negative sign for any tensor, it is enough to consider the situation that $\mathcal{A}^{n} - \mathcal{B}^{m}$ is a positive definite tensor.

If we apply $t = \log{\frac{ \lambda^n_{\mathcal{A}}}{  \lambda^m_{\mathcal{B}} }}$ to Eq.~\eqref{eq2:thm:16thm}, we have $\psi(  \lambda_{\mathcal{A}}, \lambda_{\mathcal{B}}  )$ as
\begin{eqnarray}\label{eq3:thm:16thm}
\psi(  \lambda_{\mathcal{A}}, \lambda_{\mathcal{B}}  ) = \begin{cases}
 \lambda^{\frac{n(1 - \omega)}{2}}_{\mathcal{A}} \frac{  \lambda^{n \omega}_{\mathcal{A}} -  \lambda^{m \omega}_{\mathcal{B}}   }{ \lambda^n_{\mathcal{A}} -  \lambda^m_{\mathcal{B}}  }  \lambda^{\frac{m(1 - \omega)}{2}}_{\mathcal{B}} , &  \mbox{if $\lambda_{\mathcal{A}} \neq \lambda_{\mathcal{B}}$;}\\
0, & \mbox{otherwise.}
\end{cases}
\end{eqnarray}

If we set the tensor $\mathcal{H}_t$ as
\begin{eqnarray}\label{eq4:thm:16thm}
\mathcal{H}_t = \mathcal{B}^m + t (  \mathcal{A}^n   -   \mathcal{B}^m), 
\end{eqnarray} 
then, we have 
\begin{eqnarray}\label{eq5:thm:16thm}
\lim\limits_{\delta t \rightarrow 0} \mathcal{H}^{\frac{\omega - 1}{2}}_{t + \delta t} \frac{ \mathcal{H}_{t +\delta t }    -     \mathcal{H}_{t}   }{\delta t }
\mathcal{H}^{\frac{\omega - 1}{2}}_{t} &=&  \mathcal{H}^{\frac{\omega - 1}{2}}_{t} \left( \mathcal{H}_{1}    -     \mathcal{H}_{0} \right)
\mathcal{H}^{\frac{\omega - 1}{2}}_{t}.
\end{eqnarray}

From Lemma~\ref{lma:derivative after T psi action}, and Eqs.~\eqref{eq3:thm:16thm} and~\eqref{eq5:thm:16thm}, we also have
\begin{eqnarray}\label{eq6:thm:16thm}
\frac{d}{dt}\left( \mathcal{H}_t^{\omega} \right) &=& \lim\limits_{\delta t \rightarrow 0}\frac{ \mathcal{H}_{t+ \delta t }^{\omega} - \mathcal{H}_t^{\omega}   }{ \delta t } = T_{\psi}\left(  \mathcal{H}^{\frac{\omega - 1}{2}}_{t} \left( \mathcal{H}_{1}    -     \mathcal{H}_{0} \right)
\mathcal{H}^{\frac{\omega - 1}{2}}_{t} \right).
\end{eqnarray}

By applying Lemma~\ref{lma:T psi by psi} and Corollary~\ref{cor:bound for psi norm Fourier} to the function $g(t)$ provided by Eq.~\eqref{eq2:thm:16thm}, we have
\begin{eqnarray}\label{eq7:thm:16thm}
\left\Vert \frac{d}{dt}\left( \mathcal{H}_t^{\omega} \right) \right\Vert
\leq \left( \mathbb{I}_1^N \right)^2  \left[\min\limits_{c>0} \left( \sqrt{2c}\left\Vert g(t) \right\Vert_2 + \sqrt{2/c}\left\Vert g'(t) \right\Vert_2    \right) \right]    \left\Vert   \mathcal{H}^{\frac{\omega - 1}{2}}_{t} \left( \mathcal{H}_{1}    -     \mathcal{H}_{0} \right) \mathcal{H}^{\frac{\omega - 1}{2}}_{t}  \right\Vert.
\end{eqnarray}

Because $\mathcal{H}_t - t \left(\mathcal{H}_1 - \mathcal{H}_0 \right)$ is a positive definite tensor and monotonicity of the function $t^{(1 - \omega)}$, we have 
\begin{eqnarray}\label{eq8:thm:16thm}
\left\Vert \mathcal{H}_t^{(\omega - 1)/2} \star_N  \left(   \mathcal{H}_1 - \mathcal{H}_0  \right)^{(1 - \omega)}   \star_N \mathcal{H}_t^{(\omega - 1)/2} \right\Vert \leq t^{(\omega - 1)}.
\end{eqnarray}
From Eq.~\eqref{eq7:thm:16thm} and Eq.~\eqref{eq8:thm:16thm}, we have
\begin{eqnarray}\label{eq9:thm:16thm}
\left\Vert \frac{d}{dt}\left( \mathcal{H}_t^{\omega} \right) \right\Vert &\leq& 
 \left( \mathbb{I}_1^N \right)^2  \left[\min\limits_{c>0} \left( \sqrt{2c}\left\Vert g(t) \right\Vert_2 + \sqrt{2/c}\left\Vert g'(t) \right\Vert_2    \right) \right]  \nonumber \\
&   & \times \left\Vert \mathcal{H}_t^{(\omega - 1)/2} \star_N  \left(   \mathcal{H}_1 - \mathcal{H}_0  \right)  \star_N \mathcal{H}_t^{(\omega - 1)/2} \right\Vert \nonumber \\
&\leq& \left( \mathbb{I}_1^N \right)^2  \left[\min\limits_{c>0} \left( \sqrt{2c}\left\Vert g(t) \right\Vert_2 + \sqrt{2/c}\left\Vert g'(t) \right\Vert_2    \right) \right]  \nonumber \\
&   & \times t^{(\omega - 1)} \left\Vert   \mathcal{H}_1 - \mathcal{H}_0  \right\Vert^{\omega}.
\end{eqnarray}
Therefore, we have
\begin{eqnarray}\label{eq10:thm:16thm}
\left\Vert \mathcal{A}^{n \omega} -   \mathcal{B}^{m \omega} \right\Vert &=& \left\Vert 
\int_0^{1}   \frac{d}{dt}\left( \mathcal{H}_t^{\omega} \right)   \right\Vert \nonumber \\
& \leq & \left( \mathbb{I}_1^N \right)^2  \left[\min\limits_{c>0} \left( \sqrt{2c}\left\Vert g(t) \right\Vert_2 + \sqrt{2/c}\left\Vert g'(t) \right\Vert_2    \right) \right]  \left\Vert   \mathcal{H}_1 - \mathcal{H}_0  \right\Vert^{\omega}  \int_0^1 t^{\omega - 1} dt \nonumber \\
&=_1& \frac{\left( \mathbb{I}_1^N \right)^2}{\omega}  \left[\min\limits_{c>0} \left( \sqrt{2c}\left\Vert g(t) \right\Vert_2 + \sqrt{2/c}\left\Vert g'(t) \right\Vert_2    \right) \right] \left\Vert  \mathcal{A}^{n} -   \mathcal{B}^{m}  \right\Vert^{\omega},
\end{eqnarray}
where $=_1$ is obtained by using $\mathcal{H}_1 = \mathcal{A}^n$ and $\mathcal{H}_0 = \mathcal{B}^m$  from Eq.~\eqref{eq4:thm:16thm}. This theorem is proved by applying Markov inequality to Eq.~\eqref{eq10:thm:16thm}.
$\hfill \Box$

If $m=n=1$, Theorem~\ref{thm:16thm} becomes the tail bound for Birman-Kopilenko-Solomyak inequality~\cite{birman1975estimates}. 

Following Theorem~\ref{thm:17lma} will be another tail bound for new random tensors inequality based on PDTI. 
\begin{theorem}\label{thm:17lma}
Let $\mathcal{A}, \mathcal{B} \in \mathbb{C}^{I_1 \times \cdots \times I_N \times I_1 \times \cdots \times I_N} $ be random Hermitian tensors and $\mathcal{X} \in \mathbb{C}^{I_1 \times \cdots \times I_N \times I_1 \times \cdots \times I_N} $ be a Hermitian tensor. For two real numbers $\alpha, \beta$ such that $0 \leq \alpha, \beta \leq 1$, and two positive integers $m,n$, we have
\begin{eqnarray}\label{eq1:thm:17lma}
\mathrm{Pr}\left( \left\Vert \mathcal{A}^{ \frac{m(1 + \alpha)}{2} }\star_N \mathcal{X} \star_N   \mathcal{B}^{ \frac{n(1 - \alpha)}{2} } -  \mathcal{A}^{  \frac{m(1 - \beta)}{2}   } \star_N \mathcal{X} \star_N   \mathcal{B}^{ \frac{n(1 + \beta)}{2}}  \right\Vert \geq \theta \right) ~~~~~~~ \leq \nonumber \\  
\frac{    \left( \mathbb{I}_1^N \right)^2 }{\theta}  \left[\min\limits_{c>0} \left( \sqrt{2c}\left\Vert g(t) \right\Vert_2 + \sqrt{2/c}\left\Vert g'(t) \right\Vert_2    \right) \right]   \mathbb{E}\left( \left\Vert \mathcal{A}^m \star_N \mathcal{X} - \mathcal{X} \star_N \mathcal{B}^n \right\Vert \right)
 \end{eqnarray}
where $g(t)$ is
\begin{eqnarray}\label{eq2:thm:17lma}
g(t) \define  \frac{\exp\left(\frac{\alpha t}{2}\right) -  \exp\left(\frac{-\beta t}{2}\right)    }{\exp\left( \frac{t}{2} \right) -  \exp\left(\frac{-t}{2} \right)   }
\end{eqnarray}
\end{theorem}
\textbf{Proof:}

From Theorem~\ref{thm:perturbation theory}, we have 
\begin{eqnarray}\label{eq3:thm:17lma}
\mathcal{A}^{   \frac{m(1 + \alpha)}{2}   } \star_N \mathcal{X} \star_N   \mathcal{B}^{\frac{n(1 - \alpha)}{2}} -  \mathcal{A}^{   \frac{m(1 - \beta)}{2}    } \star_N \mathcal{X} \star_N   \mathcal{B}^{  \frac{n(1 + \beta)}{2}    } =T_{\psi}\left(  \mathcal{A}^m \star_N \mathcal{X} - \mathcal{X} \star_N \mathcal{B}^n\right),
\end{eqnarray} 
where $\psi$ is 
\begin{eqnarray}\label{eq4:thm:17lma}
\psi(\lambda_{\mathcal{A}}, \lambda_{\mathcal{B}}) =  \lambda^{\frac{m(1 - \beta)}{2}}_{\mathcal{A}}  \frac{  \lambda^{\frac{m(\alpha + \beta)}{2}}_{\mathcal{A}} - \lambda^{\frac{n(\alpha + \beta)}{2}}_{\mathcal{B}}  }{ \lambda^m_{\mathcal{A}}  - \lambda^n_{\mathcal{B}} }  \lambda^{\frac{n(1 - \alpha)}{2}}_{\mathcal{B}} .
\end{eqnarray}
Then, Eq.~\eqref{eq4:thm:17lma} will be obtained by setting $t = \log \frac{\lambda^m_{\mathcal{A}}}{ \lambda^n_{\mathcal{B}}}$ in Eq.~\eqref{eq2:thm:17lma}.

By applying Lemma~\ref{lma:T psi by psi} and Corollary~\ref{cor:bound for psi norm Fourier} to the function $g(t)$ provided by Eq.~\eqref{eq2:thm:17lma}, we have
\begin{eqnarray}
\left\Vert \mathcal{A}^{ \frac{m(1 + \alpha)}{2} }\star_N \mathcal{X} \star_N   \mathcal{B}^{ \frac{n(1 - \alpha)}{2} } -  \mathcal{A}^{  \frac{m(1 - \beta)}{2}   } \star_N \mathcal{X} \star_N   \mathcal{B}^{ \frac{n(1 + \beta)}{2}}  \right\Vert  \leq 
~~~~~~~~~~~~~~~~~~~~ \nonumber \\
\left( \mathbb{I}_1^N \right)^2  \left[\min\limits_{c>0} \left( \sqrt{2c}\left\Vert g(t) \right\Vert_2 + \sqrt{2/c}\left\Vert g'(t) \right\Vert_2    \right) \right]    \left\Vert \mathcal{A}^m \star_N \mathcal{X} - \mathcal{X} \star_N \mathcal{B}^n\right\Vert.
\end{eqnarray}
Therefore, we have
\begin{eqnarray}\label{eq5:thm:17lma}
\lefteqn{\mathrm{Pr}\left( \left\Vert \mathcal{A}^{ \frac{m(1 + \alpha)}{2} }\star_N \mathcal{X} \star_N   \mathcal{B}^{ \frac{n(1 - \alpha)}{2} } -  \mathcal{A}^{  \frac{m(1 - \beta)}{2}   } \star_N \mathcal{X} \star_N   \mathcal{B}^{ \frac{n(1 + \beta)}{2}}  \right\Vert \geq \theta \right) } \nonumber \\
&\leq&  
\mathrm{Pr}\left( \left\{ \left( \mathbb{I}_1^N \right)^2  \left[\min\limits_{c>0} \left( \sqrt{2c}\left\Vert g(t) \right\Vert_2 + \sqrt{2/c}\left\Vert g'(t) \right\Vert_2    \right) \right]   \left\Vert \mathcal{A}^m \star_N \mathcal{X} - \mathcal{X} \star_N \mathcal{B}^n \right\Vert \right\} \geq \theta \right) \nonumber \\
&=& \mathrm{Pr}\left(    \left\Vert \mathcal{A}^m \star_N \mathcal{X} - \mathcal{X} \star_N \mathcal{B}^n \right\Vert  \geq \frac{\theta}{    \left( \mathbb{I}_1^N \right)^2  \left[\min\limits_{c>0} \left( \sqrt{2c}\left\Vert g(t) \right\Vert_2 + \sqrt{2/c}\left\Vert g'(t) \right\Vert_2    \right) \right]       }     \right).
\end{eqnarray}
This theorem is proved by applying Markov inequality to Eq.~\eqref{eq5:thm:17lma}.
$\hfill \Box$

Following corollary is obtained by applying Corollary~\ref{cor:perturbation theory} to the same conditions of Theorem~\ref{thm:11lma} for the tensor $\mathcal{A}^{ \frac{m(1 + \alpha)}{2} }\star_N \mathcal{X} \star_N   \mathcal{B}^{ \frac{n(1 - \alpha)}{2} } + \mathcal{A}^{  \frac{m(1 - \beta)}{2}   } \star_N \mathcal{X} \star_N   \mathcal{B}^{ \frac{n(1 + \beta)}{2}} $. We will skip the proof here due to the similarity of the proof provided by Theorem~\ref{thm:17lma}.

\begin{corollary}\label{cor:17lma}
Let $\mathcal{A}, \mathcal{B} \in \mathbb{C}^{I_1 \times \cdots \times I_N \times I_1 \times \cdots \times I_N} $ be random Hermitian tensors and $\mathcal{X} \in \mathbb{C}^{I_1 \times \cdots \times I_N \times I_1 \times \cdots \times I_N} $ be a Hermitian tensor. For two real numbers $\alpha, \beta$ such that $0 \leq \alpha, \beta \leq 1$, and two positive integers $m,n$, we have
\begin{eqnarray}\label{eq1:cor:17lma}
\mathrm{Pr}\left( \left\Vert  \mathcal{A}^{ \frac{m(1 + \alpha)}{2} }\star_N \mathcal{X} \star_N   \mathcal{B}^{ \frac{n(1 - \alpha)}{2} } +  \mathcal{A}^{  \frac{m(1 - \beta)}{2}   } \star_N \mathcal{X} \star_N   \mathcal{B}^{ \frac{n(1 + \beta)}{2}}  \right\Vert \geq \theta \right) \leq \nonumber  \\
\frac{    \left( \mathbb{I}_1^N \right)^2 }{\theta}  \left[\min\limits_{c>0} \left( \sqrt{2c}\left\Vert g(t) \right\Vert_2 + \sqrt{2/c}\left\Vert g'(t) \right\Vert_2    \right) \right] \mathbb{E}\left( \left\Vert \mathcal{A}^m \star_N \mathcal{X} + \mathcal{X} \star_N \mathcal{B}^n \right\Vert \right),
 \end{eqnarray}
where $g(t)$ is
\begin{eqnarray}\label{eq2:cor:17lma}
g(t) \define  \frac{\exp\left(\frac{\alpha t}{2}\right) +\exp\left(\frac{-\beta t}{2}\right)    }{\exp\left( \frac{t}{2} \right) + \exp\left(\frac{-t}{2} \right)   }.
\end{eqnarray}
\end{corollary}

\section{Conclusions}\label{sec:Conclusions}

In this work, we extended our previous work about DTI to PDTI by deriving the upper bound for PDTI norm and new perturbation formula for a PDTI tensor. We also studied the convergence property of random PDTI and applied this property to characterize the tensor variation after the action of PDTI. With these new instruments, we are able to build new tail bounds for random tensors. We believe the proposed random PDTI and related tools can be applied to other fields of mathematics, e.g., noncommutative geometry.  

\bibliographystyle{IEEETran}
\bibliography{paraDTI_Bib}

\end{document}